\documentclass[preprint]{elsarticle}
\usepackage[nolists]{endfloat}

\usepackage[top=2.5cm,bottom=2.5cm,left=2cm,right=2cm]{geometry}
\usepackage[dvips]{color}

\usepackage{hyperref}
\usepackage{graphicx,subfigure}
\usepackage{amssymb,verbatim,amsthm}
\newtheorem{theorem}{Theorem}[section]
 
\newtheorem{proposition}[theorem]{Proposition}

\newtheorem{assumption}[theorem]{Assumption}
\newtheorem{remark}[theorem]{Remark}

\usepackage{natbib}
\newcommand{\N}{\mathbb{N}}
\newcommand{\Z}{\mathbb{Z}}

\newcommand{\R}{\mathbb{R}}

\newcommand{\hU}{\widehat{U}}
\newcommand{\cL}{\mathcal{L}}

\newcommand{\E}{\mathbb{E}}
\newcommand{\T}{\mathbb{T}}

\newcommand{\eps}{\sigma}
\newcommand{\Laplace}{\Delta}

\newcommand{\Pl}{P_{\lambda}}
\newcommand{\Ql}{Q_{\lambda}}
\newcommand{\Bl}{B_0({\eta})}
\newcommand{\Wl}{W_{\lambda}}
\newcommand{\Wc}{W_{\lambda}^c}

\newcommand{\cA}{\mathcal{A}}
\newcommand{\cN}{\mathcal{N}}
\newcommand{\cB}{\mathcal{B}}
\newcommand{\cC}{\mathcal{C}}
\newcommand{\ls}{\lambda^{\star}}
\newcommand{\Be}{B_j}
\newcommand{\hu}{\widehat {u}}
\newcommand{\hv}{\widehat {v}}
\newcommand{\hw}{\widehat {w}}

\newcommand{\rd}{\mathrm{d}}
\newcommand{\pd}{\partial}
\newcommand{\abs}[1]{| #1 |}
\newcommand{\eref}[1]{Eq. ( \ref{#1} )}

\newcommand{\norm}[1]{\| #1 \|}

\newcommand{\inner}[2]{\langle #1 , #2 \rangle}

\begin{document}

\title{Accuracy and Stability of Filters for
Dissipative PDEs}
%\affiliation{Warwick Mathematics Institute, University of Warwick, Coventry CV4 7AL, UK}}
\author{C.\ E.\ A.\ Brett, %}
K.\ F.\ Lam, K.\ J.\ H.\ Law, % \affil{1} \corrauth 
D.\ S.\ McCormick, M.\ R.\ Scott 
and A.\ M.\ Stuart}
% \affil{1}}
%\address{\affilnum{1}
\address{
Warwick Mathematics Institute, University of Warwick,
  Coventry CV4 7AL, UK}

%\address{$^1$ Mathematics Faculty, Open University, Milton Keynes MK7~6AA, UK} \address{$^2$ Department% of Mathematics, Imperial College, Prince Consort Road, London SW7~2BZ, UK} \address{$^3$ Department of C%omputer Science,
%University College London, Gower Street, London WC1E~6BT, UK} 
\ead{k.j.h.law@warwick.ac.uk, a.m.stuart@warwick.ac.uk}

\begin{abstract}
Data assimilation methodologies are designed 
to incorporate noisy observations
of a physical system into an underlying model 
in order to infer the properties of the state 
of the system. 
Filters refer to a class of data assimilation
algorithms designed to update the estimation of the state
in an on-line fashion, as data is acquired sequentially.
For linear problems subject
to Gaussian noise, filtering can be performed exactly using
the Kalman filter. For nonlinear systems filtering can be
approximated in a systematic way by particle filters.
However in high dimensions these particle filtering methods
can break down. Hence, for the large nonlinear systems 
arising in applications such as oceanography
and weather forecasting, 
various {\em ad hoc} filters are used,
mostly based on making Gaussian approximations. The purpose
of this work is to study the accuracy and
stability properties of these 
{\em ad hoc} filters. We work in the context of the 2D 
incompressible Navier-Stokes equation, although the
ideas readily generalize to a range of dissipative 
partial differential equations (PDEs).
By working in this infinite dimensional setting we provide an
analysis which is useful for the understanding
of high dimensional filtering, and is robust to mesh-refinement.
We describe theoretical results showing that, in the
small observational noise limit, the filters
can be tuned to perform accurately in tracking the signal
itself (filter accuracy), 
provided the system is observed in a sufficiently large
low dimensional space; roughly speaking this space
should be large enough to contain the unstable modes of the
linearized dynamics. The tuning corresponds to what is known
as {\em variance inflation} in the applied literature.
Numerical results are given which
illustrate the theory. The positive results herein concerning
filter stability complement recent numerical studies which 
demonstrate that the {\em ad hoc} filters can perform poorly in 
reproducing statistical variation about the true signal.
\end{abstract}

%\pacs{}%1315, 9440T} \submitto{\JPG} 
\date{\today}

\maketitle

\section{Introduction}
\label{sec:intro}

Assimilating large  
data sets into mathematical models of time-evolving
systems presents a major challenge in a wide range of 
applications. Since the data and the model are often 
uncertain, a natural overarching framework
for the formulation of such problems is that of Bayesian
statistics. However, for high dimensional models,
investigation of the Bayesian posterior distribution
of model state given data is not computationally feasible
in on-line situations. For this reason various {\em ad hoc} 
filters are used. The purpose of this paper is to provide an
analysis of such filters. 
%fully Bayesian approaches are not computationally feasible
%in on-line situations and various {\em ad hoc} filters
%are used. The purpose of this paper is to provide an
%analysis of such filters. 

The paradigmatic example of data assimilation is weather 
forecasting: computational models to predict the state 
of the atmosphere currently involving 
on the order of ${\mathcal O}(10^{8})$ 
unknowns but these models must be run with an initial 
condition which is only known incompletely. This is compensated
for by a large number, currently on the order of 
${\mathcal O}(10^{6})$, partial observations of the
atmosphere at subsequent times. Filters are widely used to
make forecasts which combine the mathematical model
of the atmosphere and the data to make predictions. 
Indeed the particular method of data assimilation which we 
study here includes, as a special case, the algorithm 
commonly known as 3DVAR. This method originated in weather 
forecasting. It was first proposed %by Lorenc 
at the UK 
Met Office in 1986 \cite{article:Lorenc1986}, and was developed 
by the US National Oceanic and Atmospheric Administration 
(NOAA) soon thereafter; see \cite{article:Parrish1992}.
More details of the implementation of 3DVAR by the UK Met 
Office can be found in \cite{article:Lorenc2000}, and by 
the European Centre for Medium-Range Weather Forecasts (ECMWF) 
in \cite{article:Courtier1998}. The 3DVAR algorithm is 
prototypical of the many more sophisticated filters which are 
now widely used in practice  and 
it is thus natural to study it. The reader should be aware,
however, that the development of new filters is a very
active area and that the analysis here constitutes an initial
step towards the analyses required for these  more
recently developed algorithms. For insight into some
of these more sophisticated filters see 
\cite{toth1997ensemble,evensen2009data,VL09,
harlim2008filtering,majda2010mathematical} and the references therein.

Filtering can be performed exactly
for linear systems subject to Gaussian noise: 
the Kalman filter \cite{harvey1991forecasting}.
For nonlinear or non-Gaussian scenarios the particle
filter \cite{doucet2001sequential}
may be used and provably approximates the desired
probability distribution as the number of particles
is increased \cite{bain2008fundamentals}. However
in practice this method performs poorly in high
dimensional systems \cite{Bickel}. Whilst there
is considerable research activity
aimed at overcoming this degeneration 
\cite{van2010nonlinear, chorin2010implicit},
the methodology cannot yet be viewed
as a provably accurate tool within the context of the
high dimensional problems arising in geophysical
data assimilation.  In order to circumvent problems associated
with the representation of high dimensional probability
distributions some form of {\em ad hoc}
Gaussian approximation is typically
used to create practical filters, and the 3DVAR method
which we analyze here is perhaps the simplest example of this.
This {\em ad hoc} filters may also be viewed
in the framework of nonlinear control theory.
Proving filter stability and accuracy has a long history
in this field and the paper \cite{tarn76} is 
closely related to the work we develop here. However
we work in an infinite dimensional setting, in order
to directly confront the high dimensionality of
many current real-world filtering applications,
and this brings new issues to the problem of
establishing filter stability and accuracy; 
overcoming these problems provides the focus of this paper.

In the paper \cite{lawstuart} a wide range of Gaussian
approximate filters, including 3DVAR,
are evaluated by comparing the
distributions they produce with a highly accurate
(and impractical in realistic online scenarios) 
MCMC simulation of the
desired distributions. The conclusion of that work
is that the Gaussian approximate filters perform well
in tracking the mean of the desired distribution, but poorly
in terms of statistical variability about the mean. In
this paper we provide a theoretical analysis of the
ability of the filters to estimate the mean state accurately.
Although filtering is widely used in practice, much of the 
analysis of it, in the context of fluid mechanics, works with 
finite-dimensional dynamical models. Our aim is to work directly
with a PDE model relevant in fluid mechanics, the Navier-Stokes 
equation, and thereby confront the high-dimensional nature of 
the problem head-on. 
Study of the stability of filters for data assimilation has
been a developing research area over the last few years and
the paper \cite{carrassi2008data} contains a finite
dimensional  theoretical result, 
numerical experiments in 
a variety of finite and (discretized) infinite dimensional
systems not covered by the theory, and references to
relevant applied literature.  The paper \cite{trevisan2011chaos} 
gives a review of many important aspects relating to assimilating 
the evolving unstable directions of the underlying dynamical system.
Our analysis will build in a concrete fashion
on the approach in 
\cite{olson2003determining} and \cite{hayden2011discrete}
which were amongst the first to study data assimilation directly
through PDE analysis, using ideas from the theory of determining
modes in infinite dimensional dynamical systems.
However, in contrast to those papers,
we will allow for noisy observations in our analysis. 
Nonetheless the estimates in \cite{hayden2011discrete} form
an important component of our analysis.
Furthermore the large time asymptotic results
in \cite{hayden2011discrete} constitute a limiting case
of our theory, where there is no observational noise.

The presentation will be organized as follows: in Section
\ref{sec:combine} we introduce the Navier-Stokes equation 
as the forward model of interest and formulate the inverse
problem of estimating the velocity field given partial
noisy observations. This leads to a family of
filters for the velocity field which have the
form of a non-autonomous dynamical system which blends
the forward model with the data in a judicious fashion;
Theorem \ref{t:min} describes this dynamical system
via the solution of a sequence of inverse problems. 
In Section \ref{sec:stability} we introduce notions of
stability and prove the Main Theorems \ref{t:m}
and \ref{t:mz} concerning filter stability 
and accuracy for sufficiently small observational noise. 
%section \ref{sec:stability} also contains derivation
%of the stochastic PDE (SPDE) and deterministic PDE
%which may be used to
%study filter stability for 3DVAR
%when data is aquired frequently
%in time  and the observational noise is large;
In Section \ref{sec:numerics}
we present numerical results which corroborate the analysis;
and finally, in Section \ref{sec:conclusions} we present conclusions.

\section{Combining Model and Data}
\label{sec:combine}

In  subsection \ref{ssec:NSE}
we describe the forward model that we employ
throughout the remainder of the paper: the Navier Stokes
equation on a two dimensional torus.
Then, in subsection \ref{ssec:filters}, we
describe the observational data model that we employ; using
this we apply Tikhonov-Phillips regularization to derive
the filter which we use to combine
model and data.  Subsection \ref{ssec:3DVar} contains 
a specific example of this filter, used later in
the paper for our numerical illustrations. This
example constitutes a particular instance of the 3DVAR method.

\subsection{Forward Model: Navier-Stokes equation}
\label{ssec:NSE}

In this section we establish a notation for,
and describe the properties of, the 
Navier-Stokes equation. This is the forward model which
underlies the inverse problem which we study in this paper.
We consider the 2D Navier-Stokes equation on the torus
$\T^{2} := [0,L) \times [0,L)$ with periodic boundary conditions:
%\begin{subequations}
\begin{eqnarray*}
\begin{array}{cccc}
%\begin{align}
\pd_{t}u - \nu \Laplace u + u \cdot \nabla u + \nabla p &=& f 
& {\rm for~ all ~} 
%\text{for all } 
(x, t) \in \T^{2} \times (0, \infty), \\
\nabla \cdot u &=& 0 &{\rm for~ all~ }
%\text{for all }
(x, t) \in \T^{2} \times (0, \infty), \\
u(x, 0) &=& u_{0}(x) &{\rm for~ all ~}
%\text{for all } 
x \in \T^{2}.
%\end{align}
\end{array}
\label{eq:NSE}
\end{eqnarray*}
%\end{subequations}

Here $u \colon \T^{2} \times (0, \infty) \to \R^{2}$ is a time-dependent vector field representing the velocity, $p \colon \T^{2} \times (0,\infty) \to \R$ is a time-dependent scalar field representing the pressure, $f \colon \T^{2} \to \R^{2}$ is a vector 
field representing the forcing (which we assume to be
time-independent for simplicity), and $\nu$ is the viscosity. 
In numerical simulations (see section~\ref{sec:numerics}), we
typically represent the solution via  
the vorticity $w$ and stream function $\zeta$; these 
are related through $u = \nabla^{\perp} \zeta$ and $w=\nabla^{\perp} \cdot u$, 
where $\nabla^{\perp} = (\pd_{2}, -\pd_{1})^{\mathrm{T}}$.
We define
%\[
$${\mathcal H}:= \left\{ {\rm {trigonometric\,polynomials\,}} 
u:\T^2 \to {\mathbb R}^2\,\Bigl|\, \nabla \cdot u = 0, \,\int_{\T^{2}} u(x) \, \rd x = 0 \right\}
$$
%\] 
and $H$ as the closure of ${\mathcal H}$ with respect to the
$(L^{2}(\T^{2}))^{2}$ norm. We define $P:(L^{2}(\T^{2}))^{2} 
\to H$ to be the Leray-Helmholtz orthogonal projector.

Given $k = (k_{1}, k_{2})^{\mathrm{T}}$, define $k^{\perp} := (k_{2}, -k_{1})^{\mathrm{T}}$. Then an orthonormal basis for
$H$ is given by $\psi_{k} \colon \R^{2} \to \R^{2}$, where 
\begin{equation}
\psi_{k} (x) := \frac{k^{\perp}}{|k|} \exp\Bigl(\frac{2 \pi i k \cdot
  x}{L}\Bigr)\quad {\rm for} \quad k \in \Z^{2} \setminus \{0\}.
\label{psik}
\end{equation} 
Thus for $u \in H$ we may write
$$u = \sum_{k \in \Z^{2} \setminus \{0\}} u_{k}(t) \psi_{k}(x)$$
where, since $u$ is a real-valued function, we have the 
reality constraint $u_{-k} = - \overline{u_{k}}.$
We then define the projection operators $\Pl: H \to H$ 
and $\Ql:H \to H$ by
$$\Pl u = \sum_{|2\pi k|^2 <\lambda L^2} u_{k}(t) \psi_{k}(x),
%$$\Pl u = \sum_{|k|^2 <\lambda/\lambda_1} u_{k}(t) \psi_{k}(x),
\quad \Ql=I-\Pl.$$
We let $\Wl=\Pl H$ and $\Wc=\Ql H.$

Using the Fourier decomposition of $u$, we define the
fractional Sobolev spaces
\begin{equation}
\label{eq:Hs} 
H^{s}:= \left\{ u \in H : \sum_{k \in \Z^{2} \setminus \{0\}} (4\pi^{2}\abs{k}^{2})^{s}\abs{u_{k}}^{2} < \infty \right\}
\end{equation}
with the norm $\norm{u}_{s}:= \bigl(\sum_{k} (4\pi^{2}\abs{k}^{2})^{s}
\abs{u_{k}}^{2}\bigr)^{1/2}$, where $s \in \R$. We use
the abbreviated notation $\norm{u}$ for the norm on $H^1$,
and $|\cdot|$ for the norm on $H=H^0$. 

Applying the projection $P$ to the Navier-Stokes
equation we may write it as an ODE (ordinary differential
equation) in $H$; see
\cite{constantin1988navier, temam1995navier, book:Robinson2001} 
for details. This ODE takes the form
\begin{equation}
\frac{\rd u}{\rd t} + \nu Au + \cB(u, u) = f, \quad u(0)=u_0.
\label{eq:nse}
\end{equation}
Here, $A = -P \Laplace$ is the Stokes operator, 
the term $\cB(u,u) = P(u \cdot \nabla u)$ is the bilinear form 
found by projecting the nonlinear term $u \cdot \nabla u$ onto 
$H$ and finally, with abuse of notation, $f$ is the original 
forcing, projected into $H$. 
We note that $A$ is diagonalized in the basis comprised
of the $\{\psi_k\}_{k \in {\Z}^2\backslash\{0\}}$, 
on $H$, and the smallest eigenvalue of $A$
is $\lambda_1=4\pi^2/L^2.$
The following proposition
is a classical result which implies the existence
of a dissipative semigroup for the ODE \eref{eq:nse}. 
See Theorems 9.5 and 12.5 in \cite{book:Robinson2001}
for a concise overview and \cite{temam1995navier,book:Temam1997}
for further details.

\begin{proposition} 
\label{prop:1}
Assume that $u_0 \in H^1$ and $f \in H$.
Then \eref{eq:nse} has a unique strong solution on $t
\in [0,T]$ for any $T>0:$
$$u \in L^{\infty}\bigl((0,T);H^1\bigr)\cap L^{2}\bigl((0,T);D(A)\bigr),\quad \frac{du}{dt} \in L^{2}\bigl((0,T);H\bigr).$$
Furthermore the equation has a global attractor $\cA$
and there is $K>0$ such that,
if $u_0 \in \cA$, then $\sup_{t \ge 0}\|u(t)\|^2 \le K.$
\end{proposition}

We let $\{\Psi(\cdot,t): H^1 \to H^1\}_{t \ge 0}$ 
denote the semigroup of 
solution operators for the equation \eref{eq:nse} through $t$
time units. We note that by
working with weak solutions,
$\Psi(\cdot,t)$ can be extended to act on larger spaces $H^s$,
with $s \in [0,1)$, under the same assumption on $f$;
see Theorem 9.4 in \cite{book:Robinson2001}.  
We will, on occasion, use this extension of $\Psi(\cdot,t)$
to act on larger spaces. 
%In particular we use
%it in the following proposition which follows from
%Lemmas 5.3 and 5.5 in \cite{article:Cotter2009}.
%This proposition is key to enabling us to show that the 
%Bayesian inverse problem, which underlies the
%filtering problem of interest, is well-posed.

%\begin{proposition} 
%\label{prop:2}
%Let $s \in (0,1]$ and $t_0>0.$ Then,
%for all $t>t_0$, and $u,v \in H$,
%\begin{align*}
%\begin{eqnarray*}
%\begin{array}{ll}
%\|\Psi(u,t)\|_{s}^2 &\le  t_0^{-s}C\bigl(|f|^2+|u|^2)\\
%\|\Psi(u,t)-\Psi(v,t)\|_{s}^2 &\le  t_0^{-s}C\bigl(|u|,|u-v|,|f|\bigr)
%|u-v|^2.
%\end{align*}
%\end{array}
%\end{eqnarray*}
%\end{proposition}

The key properties of the Navier-Stokes
equation that drive our analysis of the filters are summarized
in the following proposition, taken from the paper
\cite{hayden2011discrete}.
To this end, define $\Psi(\cdot)=\Psi(\cdot,h)$
for some fixed $h>0.$ 
Note that the statement here is closely related
to the {\em squeezing property} \cite{book:Robinson2001}
of the Navier-Stokes equation, a property 
employed in a wide range of applied contexts.
Furthermore, many other dissipative PDEs are
known to satisfy similar properties.

\begin{proposition} 
\label{prop:3}
Let $u \in \cA$ and $v \in H^1$. There is
$\beta=\beta(|f|,L,\nu)>0$ such that
\begin{equation}
\label{eq:1}
\|\Psi(u)-\Psi(v)\|^2 \le \exp(\beta h)\|u-v\|^2.
\end{equation}
Now let $\|u-v\| \le R$ and assume that
$$\lambda>\ls:=\frac{9c^{8/3}}{\lambda_1^{\frac13}}\Bigl(
\frac{2K^\frac12+R^{\frac12}}{\nu}\Bigr)^{8/3}$$
where $c$ is a dimensionless positive constant and
$K$ is the constant appearing in 
Proposition \ref{prop:1} above.
Then there exists $t^*=t^*(|f|,L,\nu,\lambda,R)$
with the property that, for all $h \in (0,t^*],$ 
there is $\gamma \in (0,1)$ such that
\begin{equation}
\label{eq:2}
\|Q_{\lambda}\bigl(\Psi(u)-\Psi(v)\bigr)\|^2 \le \gamma^2\|u-v\|^2.
\end{equation}
\end{proposition}

\begin{proof}
The first statement is simply Theorem 3.8 from 
\cite{hayden2011discrete}. The second statement
follows from the proof of Theorem 3.9 in the same paper, modified at the
end to reflect the fact that, in our setting, 
$\Pl \delta(0) \ne 0.$
Note also that the constant $\lambda$ appearing on the right
hand side of the lower bound for $\lambda$ in the 
statement of Theorem 3.9 in \cite{hayden2011discrete}
should be $\lambda_1$ (as the
proof that follows in \cite{hayden2011discrete} shows) 
and that use of definition of $K$ 
(see Theorem 3.6 of that paper) 
allows rewrite in terms of $K$ -- indeed
the proof in that paper is expressed in terms of $K$.
%proof that follows in their paper shows) and that use of
%their definition of $K$ (see their Theorem 3.6) allows rewrite in terms of $K$ -- indeed
%their proof is expressed in terms of $K$.
\end{proof}

\subsection{Inverse Problem: Filtering}
\label{ssec:filters}

In this section we describe the basic problem of {\em filtering}
the Navier Stokes equation \eref{eq:nse}: estimating
properties of the state of the system sequentially
from partial, noisy, sequential observations of the state. 
Throughout the following we write $\N$ for the natural 
numbers $\{ 1, 2, 3, \dots \}$, and $\Z^{+} := \N \cup \{ 0 \}$ for the non-negative integers $\{ 0, 1, 2, 3, \dots \}$.
Let $H$ be the Hilbert spaces with inner product 
$\inner{\cdot}{\cdot}$ and norm $|\cdot|$ defined
via \eref{eq:Hs} with $s=0$. 
If $\cL^{-1}$ is a self-adjoint, positive-definite compact operator on
$H$, which therefore has a symmetric square root $\cL^{-1/2}$,
we define
\[
\inner{\cdot}{\cdot}_{\cL} := \inner{\cdot}{\cL^{-1} \cdot}, \qquad \norm{\cdot}_{\cL} := |\cL^{-1/2} \cdot|.
\]
Recall that we have defined $\Psi(\cdot)=\Psi(\cdot,h)$
for some fixed $h>0.$ 
We let $X$ denote $H^1$ and
define $\{ u_j \}_{j \in \Z^{+}}$, $u_j \in X$ by 
\footnote{With abuse of notation, subscripts $j$ will indicate
times, while subscripts $k$ will denote Fourier coefficients as before
in order to avoid confusion. The meaning should also be clear in context.}
\begin{equation}
\label{eqn:ModelWithoutNoise}
u_{j+1} := \Psi(u_j).
\end{equation}
This discrete dynamical system is well-defined by
virtue of Proposition \ref{prop:3}.
Thus $u_j=u(jh)$ where $u$ is the solution of \eref{eq:nse}.

Our interest is in determining $u_j$ from noisy
observations of $\Pl u_j$.
We now let $\{ {\xi}_j \}_{j \in \N}$ be a noise sequence in $\Wl$ 
which perturbs the sequence $\{\Pl u_j \}_{j \in \N}$
to generate the observation sequence $\{ y_j \}_{j \in \N}$ 
in $\Wl$ given by
\begin{equation}
\label{eqn:Observation}
y_{j+1} := \Pl u_{j+1} + {\xi}_{j+1}, \quad j \in \Z^{+}.
\end{equation}
This observation operator allows us to study
partially observed infinite dimensional systems in
a clean fashion and the resulting analysis will be
a useful building block for the study of other
partial observations such 
as pointwise, or smoothed, observations on a regular grid 
in physical space. 

We let $Y_j=\{y_i\}_{i=1}^j$, the accumulated data 
up to time $t=jh.$
We assume that $u_0$ is not known exactly. The {\em goal}
of filtering is to determine the state
$u_j$ from the data $Y_j.$ 
The approach to the problem of blending data and model 
that we study here is based on Tikhonov-Phillips regularization. 
We find a point which represents
the best compromise between information given by the model
and by the data. More precisely, we use the model to provide
a regularization of a least squares problem designed to match
the data. This will result in a sequence $\{\hu_j\}_{j \in \Z^+}$ 
which approximates the true signal $\{u_j\}_{j \in \Z^+}$
giving rise to the data. To define the sequence 
$\{\hu_j\}_{j \in \Z^+}$ we introduce two 
sequences of operators $\{\Gamma_j\}_{j \in \N}$,
and $\{\cC_j\}_{j \in \N}$ which will be used to weight
the contributions of model and data at discrete time $j.$
Then define
\begin{equation}
\label{eq:Ij}
I_j(u)=\frac12\|y_{j+1}-\Pl u\|_{\Gamma_{j+1}}^2
+\frac12\|u-\Psi(\hu_j)\|_{\cC_{j+1}}^2.
\end{equation}
Choosing a minimizer of the functional encodes the
idea of a compromise between the model output $\Psi(\hu_j)$
and the data $y_{j+1}$, to estimate the state of the
system at time $t=(j+1)h$. The operators
$\cC_{j+1}$ and $\Gamma_{j+1}$
give appropriate weights on the two sources of information. 
With this in mind we set
\begin{equation}
\hu_{j+1}={\rm arginf}_{u \in H^1} I_j(u).
\label{eq:update}
\end{equation}
The following theorem shows that this is justified:

\begin{theorem}
\label{t:min}
Assume that $\cC_{j+1}$ and $\Gamma_{j+1}$ are positive-definite
self-adjoint operators on $H$ and $\Pl H$ respectively,
that $\cC_{j+1}$ is a bounded operator from
$H^1$ into itself and that the norm $\|\cdot\|_{\cC_{j+1}}$ is
equivalent to the $H^r$ norm for some $r \ge 1.$ Then,
if $\hu_0 \in H^1$, the iteration (\ref{eq:Ij},
\ref{eq:update}) determines a unique sequence 
$\{\hu_j\}_{j \in \N}$ in $H^1$ given by 
\begin{equation}
\label{eq:min}
\hu_{j+1}=\bigl(I-K_{j+1}\bigr)\Psi(\hu_j)+K_{j+1}y_{j+1}
\end{equation}
where
\begin{equation}
\label{eq:gain}
K_{j+1}=\cC_{j+1}\Pl^*\Bigl(\Gamma_{j+1}+\Pl\cC_{j+1}\Pl^*\Bigr)^{-1}
\Pl.
\end{equation}
%\begin{equation}
%\label{eq:min}
%\bigl(\cC_{j+1}^{-1}+\Pl^*\Gamma_j^{-1}\Pl\bigr)\hu_{j+1}=
%\cC_{j+1}^{-1}\Psi(\hu_j)+\Pl^*\Gamma_j^{-1}y_{j+1}.
%\end{equation}
\end{theorem}

\begin{proof} The proof follows by a simple induction. Assume
that $\hu_j \in H^1$. Then let $u=\Psi(\hu_j)+v$ and note that
the functional $J_j: H^{r} \to \R^+$ given by
\begin{equation}
\label{eq:Jj}
J_j(v)=\frac12\|y_{j+1}-\Pl\bigl(\Psi(\hu_j)+v\bigr)\|_{\Gamma_{j+1}}^2
+\frac12\|v\|_{\cC_{j+1}}^2
\end{equation}
is convex and weakly lower semi-continuous and hence
has a unique minimizer $\hv \in H^{r}$. 
Thus $\hu_{j+1}=\hv+\hu_{j}$ is in $H^1$ since
$\hu_j \in H^1$ and $\hv \in H^{r}$ with $r \ge 1$.
Equations (\ref{eq:min}, \ref{eq:gain}) for the minimizer 
follow from standard theory of calculus of variations.
% from the Euler-Lagrange equations.
\end{proof}

In the case where
$\Psi(\cdot)$ is a linear map, the matrix $K_{j+1}$ is
the {\em Kalman gain} matrix \cite{harvey1991forecasting},
composed with the observation operator $\Pl.$
We will also study the situation where
complete observations are made, obtained by taking $\lambda
\to \infty$ in the preceding analyses.
The observations are given by
\begin{equation}
\label{eqn:Observation2}
y_{j} := u_{j} + {\xi}_{j}, \quad j \in \Z^{+}
\end{equation}
where now $y_j, \xi_j \in H.$ 
The operator $\Gamma_j$ is now a positive self-adjoint
operator on $H$.
The filter that we study takes the form (\ref{eq:min}) 
with $\Pl$ replaced by the identity so that (\ref{eq:gain})
is replaced by
\begin{equation}
\label{eq:gain2}
K_{j+1}=\cC_{j+1}\Bigl(\Gamma_{j+1}+\cC_{j+1}\Bigr)^{-1}.
\end{equation}

If we define 
\begin{equation}
\label{eq:Bdef}
B_j=I-K_{j+1}
\end{equation}
then Theorem \ref{t:min} yields the key equation
\begin{equation}
\label{eq:meanu}
\hu_{j+1}=B_j\Psi(\hu_j)+(I-B_j)y_{j+1}.
\end{equation}
This equation and \eref{eq:min} demonstrate that the estimate of
the solution at time $j+1$ is found as an operator-convex
combination of the true dynamics applied to the estimate
of the solution at time $j$, and the data at time $j+1$.  We will
favor use of $B_j$ in what follows, rather than $K_{j+1}$,
as $B_j$ is the operator which controls the
stability, and hence accuracy, of the estimate.

This problem can be given a Bayesian formulation
in which the probability distribution of $u_j|Y_j$
is the primary object of interest. In practice, for
high dimensional systems arising in applications in the
atmospheric sciences, various {\em ad hoc}
Gaussian approximations are typically used to
approximate these probability distributions; the
reader interested in understanding the derivation of
filters from this probabilistic perspective is directed
to \cite{lsetal} for details; this unpublished technical
report is an expanded version of the  material contained
here, to incorporate the probabilistic perspective.
The mean of a Gaussian
posterior distribution has an equivalent formulation
as solution of a Tikhonov-Phillips regularized
least squares problem, as we have here, and this perspective on
data assimilation is adopted in \cite{pot12}. In
\cite{pot12} linear autonomous dynamical systems are studied and 
filter accuracy and stability results 
similar to ours are derived.

It is demonstrated numerically in \cite{lawstuart} that 
the Gaussian approximations are, in general, 
not good approximations.
More precisely, they fail to accurately capture
covariance information. However, the same numerical experiments
reveal that the methodology can perform well in
replicating the mean, if parameters are chosen correctly,
even if it is initially in error. Indeed this accurate
tracking of the mean is often achieved by means of
{\em variance inflation} -- increasing the model uncertainty,
here captured in the exogenously imposed $\cC_j$, 
in comparison with the 
data uncertainty, here captured in the $\Gamma_j$.
The purpose of the remainder of the paper is to
explain, and illustrate, this phenomenon by means of
analysis and numerical experiments.

\subsection{Example of a Filter: 3DVAR}
\label{ssec:3DVar}

The algorithm described in the previous section yields
the well-known 3DVAR method, discussed in the introduction,
when $\cC_j \equiv \cC$ and $\Gamma_j \equiv \Gamma$
for some fixed operators $\cC, \Gamma$.
To impose commutativity with $A$, we assume that 
the operators $\Gamma$ and $\cC$
are both fractional powers of the Stokes operator $A$, in
$\Wl$ and $H$ respectively.
We choose $A_0=\ell A$ (the parameter $\ell$ forms
a useful normalizing constant in the numerical experiments
of section \ref{sec:numerics}) and
set $\cC =\delta^2 A_0^{-2\zeta}$ in $H$ and 
$\Gamma = \sigma^2 A_0^{-2\beta}$ in $\Wl.$ Substituting 
into the update formula \eref{eq:meanu}
for $\hu_j$ and defining 
$\eta = \sigma / \delta,$ $\alpha = \zeta - \beta$,
$\Bl=(I+\eta^2 A_0^{2\alpha})^{-1}\eta^2 A_0^{2\alpha}$ in $\Wl$ 
then in (\ref{eq:meanu}) we have $B_j=B:\Wl \times \Wc \to \Wl \times \Wc$ the constant operator 
\begin{eqnarray}
B=\left(
\begin{array}{cc}
\Bl & 0\\
0 & I 
\end{array}
\right).
\label{eq:B}
\end{eqnarray}
Using this we obtain the mean update formula
\begin{equation}
\hu_{j+1} = B\Psi(\hu_j) + ( I- B) y_{j+1}.
\label{eq:up2}
\end{equation}
Notice that for $\cC, \Gamma$ given as above, 
the algorithm depends only on the three parameters
$\lambda, \alpha$ and $\eta$, once the constant
of proportionality $\ell$ in $A_0$ is set. The parameter $\lambda$
measures the size of the space in which observations
are made; for fixed wavevector $k$, the parameter $\eta$ 
is a measure of the
scale of the uncertainty in observations to uncertainty
in the model; and the sign of the parameter $\alpha$
determines whether, for fixed $\eta$ and asymptotically 
for large wavevectors,
the model is trusted more ($\alpha>0$) or less ($\alpha<0$)
than the data.  This can be seen by noting that if $\alpha>0$
then $B \psi_k \rightarrow \psi_k$, while if $\alpha<0$
then $B \psi_k \rightarrow 0$.''

In the case $\lambda=\infty$, the case of complete
observations where the whole velocity field 
is noisily observed, we again obtain \eref{eq:up2}, with 
$B=\Bl=(I+\eta^2 A_0^{2\alpha})^{-1}\eta^2 A_0^{2\alpha}$ in $H$.
The roles of $\eta$ and $\alpha$ are the same as in the
finite $\lambda$ (partial observations) case.

The discussion concerning parametric
dependence with respect to varying $\eta$ shows that,
for the example of 3DVAR introduced here, and
for both $\lambda$ finite and infinite, 
variance inflation, which refers to reducing
faith in the model in comparison with the data,
can be achieved by decreasing
the parameter $\eta.$ We will show that variance inflation 
does indeed improve the ability of the filter to track 
the signal.

\section{Accuracy and Stability}
\label{sec:stability}

In this section we develop conditions under which 
it is possible to prove
stability of the nonautonomous dynamical
system defined by the mean update equation \eref{eq:meanu}
and show that, after a sufficiently long time,
the true signal is accurately recovered. 
By stability we here mean that two filters
driven by the same noise observation will converge
towards the same estimate of the solution. By
accuracy we mean that when the noise
perturbing the observations is ${\mathcal O}(\epsilon)$, 
the filter 
will converge to an ${\mathcal O}(\epsilon)$ neighbourhood 
of the true signal, even if initially it is an ${\mathcal O}(1)$
distance from the true signal.
In subsection \ref{ssec:m1} we study the case of partial
observations; subsection \ref{ssec:m2} contains the
(easier) result for the case of complete observations.
The third subsection \ref{ssec:s2} shows how
our results can be applied to the specific
instance of the 3DVAR algorithm introduced
in subsection \ref{ssec:3DVar}, for any $\alpha \in {\mathbb R}$,
provided $\eta$, which is a measure of the ratio of uncertainty
in the data to uncertainty in the model, is sufficiently small:
this, then, is a result concerning variance inflation.

For simplicity, we will assume a ``truth'' which is on
the global attractor, as in \cite{hayden2011discrete}. This
is not necessary, but streamlines the presentation
as it gives an automatic uniform in time bound in $H^1$. 
Recall that $\|\cdot\|$ denotes the norm 
on $H^1$,
and $|\cdot|$ the norm on $H$; similarly we lift
$\|\cdot\|$ to denote the induced operator norm
on $H^1 \to H^1$.

It is useful to recall the filter in the
form (\ref{eq:meanu}):
\begin{equation}
\hu_{j+1}=B_j\Psi(\hu_j)+(I-B_j)y_{j+1}.
\end{equation}
It is also useful to consider a second filter
driven by the same data $\{y_j\}_{j \in \Z^+}$, but
possibly started at a different point:
\begin{equation}
\label{eq:meanw}
\hw_{j+1}=B_j\Psi(\hw_j)+(I-B_j)y_{j+1}.
\end{equation}

\subsection{Main Result: Partial Observations}
\label{ssec:m1}

In this case we will see that it is crucial that the 
observation space $\Wl$ is sufficiently large, i.e. 
that a sufficiently large number of modes are observed. 
This, combined with the contractivity in the high modes 
encapsulated in Proposition \ref{prop:3} 
from \cite{hayden2011discrete},
can be used to ensure stability if combined with variance
inflation.  
We study filters of the form given in \eref{eq:meanu}
and make the following assumption on the observations $\{y_j\}.$ 

\begin{assumption}
\label{a:1}
Consider a sequence $u_j=u(jh)$, where $u(t)$ is a solution
of \eref{eq:nse} lying on the global attractor $\cA$. 
Then, for some $\lambda \in (\lambda_1,\infty)$,
$$y_j=\Pl u_j+\xi_j$$
for some sequence $\xi_j$ satisfying $\sup_{j \ge 1}\|\xi_j\|
\le \epsilon.$
\end{assumption}

Note that this assumption, concerning uniform boundedness of the noise, is not
verified for the i.i.d. Gaussian case. However we do expect that a more involved
analysis would enable us to handle Gaussian noise, at the expense of proving
results in mean square, or in probability. Indeed in \cite{blsz12}
we study a continuous time limit of the set-up contained in this paper, in which
white noise forcing arises from the i.i.d. Gaussian noise; accuracy and stability 
results can then indeed be proved in mean square and in probability. However,
we believe that, for clarity, the assumption made in this paper enables us to convey
the important ideas in the most straightforward fashion.

We make the following assumption about the family $\{\Be\}$,
and assumed dependence on a parameter $\eta \in {\mathbb R}^+.$
Recall that the inverse of $\eta$ 
quantifies the amount of variance inflation. 

\begin{assumption}
\label{a:2}
The family of positive operators $\{\Be(\eta)\colon H^1 \to H^1\}_{j \ge 1}$ 
commute with $A$, satisfy $\sup_{j \ge 1}\|\Be(\eta)\| \le 1$, and 
$\sup_{j \ge 1}\|I-\Be(\eta)\| \le b$ for some $b \in {\mathbb R}^+,$
uniformly with respect to $\eta$.
Furthermore, $\bigl(I-B_j(\eta)\bigr)\Ql \equiv 0$ 
and there is, for all $\lambda>\lambda_1$, constant 
$c=c(\lambda)>0$ such that $\sup_{j \ge 1}\|\Pl \Be(\eta)\| \le c\eta^2.$
\end{assumption}

We now study the asymptotic behaviour of the filter
under these assumptions. 

\begin{theorem}
\label{t:m}
Let Assumptions \ref{a:1} and \ref{a:2} hold, choose
any $\hu_0, \hw_0 \in \mathbb{B}_{H^1}\bigl(u(0),r\bigr)$ and
let $(\lambda^*,t^*)$ be as given in Proposition \ref{prop:3}.
Assume that $\lambda>\lambda^*$.  Then for any $h \in (0,t^*]$
there is $\eta$ sufficiently 
small so that the sequences
$\{\hu_j\}_{j \ge 0}$, $\{\hw_j\}_{j \ge 0}$ given 
by \eref{eq:meanu}, \eref{eq:meanw} satisfy,
for some $a \in (0,1)$,
$$\|\hu_j-\hw_j\| \le a^j\|\hu_0-\hw_0\|$$
and
$$\|\hu_j-u_j\| \le a^j r+2b\epsilon\sum_{i=0}^{j-1}
a^i.$$ 
Hence
$$\limsup_{j \to \infty}\|\hu_j-u_j\| \le \frac{2b}{1-a}\epsilon.$$ 
\end{theorem}

\begin{proof} 
We prove the second, accuracy, result concerning $\|\hu_j-u_j\|.$
The stability result concerning $\|\hu_j-\hw_j\|$ is
proved similarly.
Assumption \ref{a:2} shows that
$y_{j+1}=\Pl\Psi(u_j)+\xi_{j+1}.$ 
Recall that in \eref{eq:meanu} $y_{j+1}$
has been extended to an element of $H$, by defining it to
be zero in $\Wc$, and we do the same with $\xi_{j+1}.$
Substituting the resulting expression for $y_{j+1}$ in
\eref{eq:meanu} we obtain
$$\hu_{j+1}=\Be \Psi(\hu_j)+(I-\Be)\Pl\Psi(u_j)+(I-\Be)\xi_{j+1}
$$
but since $(I-\Be)\Ql \equiv 0$ by assumption we have
\begin{equation}
\label{eq:need}
\hu_{j+1}=\Be \Psi(\hu_j)+(I-\Be)\Psi(u_j)+(I-\Be)\xi_{j+1}.
\end{equation}
Note also that
$$u_{j+1}=\Be \Psi(u_j)+(I-\Be)\Psi(u_j).$$
Subtracting gives the basic equation for error propagation,
namely
\begin{equation}
\hu_{j+1}-u_{j+1}=\Be\bigl(\Psi(\hu_j)-\Psi(u_j)\bigr)
+(I-\Be)\xi_{j+1}.
\label{eq:error}
\end{equation}

Since $\lambda>\lambda^{\star}$ the second
item in Proposition \ref{prop:3} holds.
Fix $a \in (\gamma,1)$ where $\gamma$ is defined in
Proposition \ref{prop:3}.  Assume,
for the purposes of induction, that 
$$\|\hu_j-u_j\| \le a^j r+2b\epsilon\sum_{i=0}^{j-1}
a^i.$$ 
Define $R=2r$ noting that
the inductive hypothesis implies that, for $\epsilon$
sufficiently small, 
$\|\hu_j-u_j\| \le r+2b(1-a)^{-1}\epsilon \le R.$
Applying $\Pl$ to \eref{eq:error} and using \eref{eq:1} gives
\begin{eqnarray*}
\begin{array}{cc}
%\begin{align*}
\|\Pl(\hu_{j+1}-u_{j+1})\| &\le
\|\Pl\Be\|\|\bigl(\Psi(\hu_j)-\Psi(u_j)\bigr)\|
+\|\Pl(I-\Be)\|\epsilon\\
&\le c(\lambda) \eta^2 \exp(\beta h/2)\|\hu_j-u_j\|+b\epsilon.
%\end{align*}
\end{array}
\end{eqnarray*}
Applying $\Ql$ to \eref{eq:error} and using \eref{eq:2}
gives\footnote{The term $b\epsilon$ on the right-hand side of
the final identity can here be set to zero because 
$(I-\Be)\Ql \equiv 0$; however in the
analogous proof of Theorem \ref{t:mz} it is present and
so we retain it for that reason.}
\begin{eqnarray*}
\begin{array}{cc}
%\begin{align*}
\|\Ql(\hu_{j+1}-u_{j+1})\| &\le
\|\Be\|\|\Ql\bigl(\Psi(\hu_j)-\Psi(u_j)\bigr)\|
+\|\Ql(I-\Be)\|\epsilon\\
&\le \gamma\|\hu_j-u_j\|+b\epsilon.
%\end{align*}
\end{array}
\end{eqnarray*}
Now note that, for any $w \in H^1$,
$\|w\| =\bigl(\|\Pl w\|^2+\|\Ql w\|^2\bigr)^{\frac12}
\le \|\Pl w\|+\|\Ql w\|.$
Thus, by adding the two previous inequalities, we find
that
$$\|\hu_{j+1}-u_{j+1}\| \le 
\bigl(c(\lambda) \eta^2 \exp(\beta h/2)+\gamma\bigr)
\|\hu_j-u_j\|+2b\epsilon.$$
Since $\gamma\in (0,1)$ and $a \in (\gamma,1)$, we may
choose $\eta$ sufficiently small so that 
$$\|\hu_{j+1}-u_{j+1}\| \le a\|\hu_j-u_j\|+2b\epsilon.$$
and the inductive hypothesis holds with $j \mapsto j+1$.
Taking $j \to \infty$ gives the desired result concerning
the limsup.
\end{proof}

\begin{remark}
\label{rem:n}

Note that the proof exploits the fact that
$B_j\Psi(\cdot)$ induces a contraction within 
a finite ball in $H^1$. This contraction
is established by means of the contractivity of $B_j$
in $\Wl$, via variance inflation,
and the squeezing property of $\Psi(\cdot)$ in $\Wc$,
for large enough observation space,
from Proposition \ref{prop:3}. 

There are two important conclusions from this theorem. The
first is that, even though the solution is only
observed in the low modes, there is sufficient
contraction in the high modes to obtain an error
in the entire estimated state which is of the same order
of magnitude as the error in the (low mode only) observations.
The second is that this phenomenon occurs even when
the initial estimate suffers from an ${\mathcal O}(1)$ error.
Of course this result also shows that, if the true solution starts in 
an $O(\epsilon)$ neighbourhood of the truth, then it remains there
for all positive time.

\label{rem:1}
\end{remark}

\subsection{Main Result: Complete Observations}
\label{ssec:m2}

Here we study filters of the form given in \eref{eq:meanu}
with observations given by \eref{eqn:Observation2}.
In this situation  the whole velocity field is
observed and so, intuitively, it should be
no harder to obtain stability than in the partially
observed case. The proof is in fact almost identical
to the case of partial observations, and so we omit
the details. We observe that, although there is no parameter
$\lambda$ in the problem statement itself, 
it is introduced in the proof: as in the previous
subsection, see Remark \ref{rem:n},
 the key to stability is to obtain
contraction in $\Wc$ using the squeezing property of
the Navier-Stokes equation, and contraction in
$\Wl$ using the properties of the filter to control
unstable modes. 

We make the following assumptions: 

\begin{assumption}
\label{a:1z}
Consider a sequence $u_j=u(jh)$ where $u(t)$ is a solution
of \eref{eq:nse} lying on the global attractor $\cA$. 
Then
$$y_j=u_j+\xi_j$$
for some sequence $\xi_j$ satisfying $\sup_{j \ge 1}\|\xi_j\|
\le \epsilon.$
\end{assumption}

\begin{assumption}
\label{a:2z}
The family of positive operators $\{\Be(\eta)\colon H^1 \to H^1\}_{j \ge 1}$ 
commute with $A$  $\sup_{j \ge 1}\|\Be(\eta)\| \le 1$, and 
$\sup_{j \ge 1}\|I-\Be(\eta)\| \le b$ for some 
$b \in {\mathbb R}^+,$ uniformly with respect to $\eta.$
Furthermore, for all $\lambda>\lambda_1$,  there is a constant 
$c=c(\lambda)>0$ such that 
$\sup_{j \ge 1}\|\Pl \Be(\eta)\| \le c\eta^2.$
\end{assumption}

We now study the asymptotic behaviour of the filter
under these assumptions. 

\begin{theorem}
\label{t:mz}
Let Assumptions \ref{a:1z} and \ref{a:2z} hold and choose
any $\hu_0, \hw_0 \in \mathbb{B}_{H^1}\bigl(u(0),r\bigr).$ 
Then for any $h \in (0,t^*]$, with $t^*$ given in 
Proposition \ref{prop:3}, there is $\eta$ sufficiently 
small so that the sequences
$\{\hu_j\}_{j \ge 0}$, $\{\hu_j\}_{j \ge 0}$ given 
by \eref{eq:meanu}, \eref{eq:meanw} satisfy,
for some $a \in (0,1)$,
$$\|\hu_j-\hw_j\| \le a^j\|\hu_0-\hw_0\|$$
and
$$\|\hu_j-u_j\| \le a^j r+2b\epsilon\sum_{i=0}^{j-1}
a^i.$$ 
Hence
$$\limsup_{j \to \infty}\|\hu_j-u_j\| \le \frac{2b}{1-a}\epsilon.$$ 
\end{theorem}

\begin{proof} The proof is nearly identical to that of
Theorem \ref{t:m}. Differences arise only because we have
not assumed that $(I-\Be)\Ql \equiv 0.$ This fact
arises in two places in Theorem \ref{t:m}.
The first is where we obtain \eref{eq:need}.
However in this case we directly obtain \eref{eq:need} 
since the whole velocity field is observed. The second
place it arises is already dealt with in the footnote
appearing in the proof of Theorem \ref{t:m}
when estimating the contraction
properties in $\Wc$; there we indicate that the proof is
already adjusted to allow for the situation required here. 
\end{proof}

\begin{remark}
\label{r:2}

If $\sup_{j \ge 1}\|B_j(\eta)\|<c\eta^2$ 
then the proof may be simplified considerably
%If $\|B_j\|<1$ then the proof may be simplified considerably
as it is not necessary to split the space into two parts,
$\Wl$ and $\Wc$. Instead the contraction of $B_j$ can
be used to control any expansion in $\Psi(\cdot)$, provided
$\eta$ is sufficiently small.
\end{remark}

\subsection{Example of Main Result: 3DVAR}
\label{ssec:s2}

We demonstrate that the 3DVAR algorithm from 
subsection \ref{ssec:3DVar} satisfies Assumptions
\ref{a:2} and \ref{a:2z} in the partially and
completely observed cases respectively, and hence
Theorems \ref{t:m} and \ref{t:mz} respectively 
may be applied to the resulting filters. 
In particular, the filters will locate the true signal,
provided $\eta$ is sufficiently small.  
Satisfaction of Assumptions \ref{a:2} and \ref{a:2z}
follows from the properties of 
$$\Bl=(I+\eta^2A_0^{2\alpha})^{-1}\eta^2 A_0^{2\alpha}, \quad
I-\Bl=(I+\eta^2A_0^{2\alpha})^{-1}.$$
Note that the eigenvalues of $\Bl$ are
$$\frac{\eta^2\bigl(4\ell\pi^2|k|^2\bigr)^{2\alpha}}{1+\eta^2\bigl(4\ell\pi^2|k|^2\bigr)^{2\alpha}},$$
if $A_0=\ell A.$
Clearly the spectral radius of $\Bl$ is less than or equal
to one on $\Wl$ or $H$, independently of the sign of
$\alpha.$  The difference is just that $|k|^2<\lambda/\lambda_1$
in the former, and $|k|$ is unbounded in the latter.

First we consider the partially observed situation.
We note that $B_j \equiv B$ and is given by \eref{eq:B}: 
\begin{eqnarray}
B=\left(
\begin{array}{cc}
(I+\eta^2A_0^{2\alpha})^{-1}\eta^2 A_0^{2\alpha} & 0\\
0 & I 
\end{array}
\right)
\label{be}
\end{eqnarray}
so that the Kalman gain-like matrix $I-B$ is given by
\begin{eqnarray}
I-B=\left(
\begin{array}{cc}
(I+\eta^2A_0^{2\alpha})^{-1} & 0\\
0 & 0 
\end{array}
\right).
\label{ibe}
\end{eqnarray}
From this it is clear that
$(I-B)\Ql \equiv 0.$ Furthermore, since the spectral
radius of $\Bl$ does not exceed one, the same is true of $B$.
Hence for the operator norms from $H^1$ into itself
we have $\|B\| \le 1.$ Similarly,
if $\alpha<0$ then $b:=\|I-B\|=1$, whilst if $\alpha \ge 0$
%if $\alpha<0$ then $\|I-B\|=b=1$, whilst if $\alpha \ge 0$
then $b=\Bigl(1+\eta^2(\ell \lambda_1)^{2\alpha}\Bigr)^{-1}<1.$ 
Thus Theorem \ref{t:m} applies.
Note that $P_{\lambda}B=B_0$ and that $B_0=O(\eta^2)$.

In the fully observed case we simply have $B_j \equiv B$
where $B=B_0(\eta)$ defined above on $H$. 
%$$B=(I+\eta^2A_0^{2\alpha})^{-1}\eta^2 A_0^{2\alpha}, \quad
%I-B=(I+\eta^2A_0^{2\alpha})^{-1}.$$ 
Again $\|B\| \le 1$ and
if $\alpha<0$ then $\|I-B\|=b=1$, whilst if $\alpha \ge 0$
then $b=\Bigl(1+\eta^2(\ell\lambda_1)^{2\alpha}\Bigr)^{-1}<1.$
Thus Theorem \ref{t:mz} applies. Note (see
Remark \ref{r:2}), that if $\alpha<0$ then
the proof of that theorem could be
simplified considerably because $\|B\|<1$ and in fact
$\sup_{j \ge 1}\|B\|<c\eta^2.$ 

\begin{remark}
\label{r:3} 

We observe that the key conclusion of Theorems \ref{t:m}
and \ref{t:mz} is
the asymptotic accuracy of the algorithm, when started at
distances of ${\mathcal O}(1)$. The asymptotic bound,
although of ${\mathcal O}(\epsilon)$, has constant
$\frac{2b}{1-a}$ which may exceed $1$
and so the bound may exceed the error obtained
by simply using the
observations to estimate the signal. Our
numerics will show, however, that in practice
the algorithm gives estimates
of the state which improve upon the observations.
\end{remark}

\section{Numerical Results}
\label{sec:numerics}

% \begin{figure*}
% \includegraphics[width=0.4\textwidth]{../figures/steadystate.png}
% \end{figure*}

In this section we describe a number of numerical
results designed to illustrate the range of filter
stability phenomena studied in the previous sections. 
We start, in subsection \ref{ssec:prob}, by
describing two useful bounds on the error committed
by filters; we will use these guides in the
subsequent numerics. Subsection \ref{ssec:results} describes
the common setup for all the subsequent
numerical results shown. Subsection \ref{ssec:complete}
describes these results in 
the case of complete observations in discrete time, whilst
Subsection \ref{ssec:partial} extends to the case of partial
observations, also in discrete time.

Our theoretical results have been derived under
Assumptions \ref{a:1} and \ref{a:1z} on the errors.
These are incompatible
with the assumption that the observational noise sequence
is Gaussian. This is because i.i.d Gaussian sequences
will not have finite supremum. However, in order to test the
robustness  of our theory we will conduct numerical
experiments with Gaussian noise sequences. 

\subsection{Useful Error Bounds}
\label{ssec:prob}

We describe two useful bounds on the error
which help to guide and evaluate the numerical
simulations. To derive these bounds we assume that
the observational noise sequence $\xi_j$ is i.i.d
with $\E \xi_j=0$ and $\E \xi_j \otimes \xi_j = \Gamma$.  
Then 
$$\E |\xi_j|^2={\rm tr} (\Gamma) = \sum_k g_k$$ 
where $\{g_k\}$ are the eigenvalues of the operator 
$\Gamma.$

\begin{itemize}

\item The lower bound is derived from 
(\ref{eq:error}). Using the assumed independence of
the sequence we see that 
\begin{equation}
\E |\hu_{j+1}-u_{j+1} |^2 \ge \E |(I-B_j)\xi_{j+1} |^2=
{\rm tr}\Bigl((I-B_j)\Gamma(I-B_j)^*\Bigr)
\label{eq:lowerbd}
\end{equation}

\item The upper bound on the filter error is found by
noting that a trivial filter is obtained by 
simply using the observation sequence as the
filter mean; this corresponds to setting
$B_j \equiv 0$ in (\ref{eq:meanu}). For this filter
we obtain 
\begin{equation}
\E |\hu_{j+1}-u_{j+1} |^2 \le \E |\xi_{j+1} |^2=
{\rm tr}\bigl(\Gamma\bigr)
\label{eq:upperbd}
\end{equation}
in the case of complete observations, and
\begin{equation}
\E |\hu_{j+1}-u_{j+1} |^2 \le \E |\xi_{j+1} |^2+|\Ql u_{j+1}|^2=
{\rm tr}\bigl(\Gamma\bigr)+|\Ql u_{j+1}|^2
\label{eq:upperbd2}
\end{equation}
in the case of incomplete observations. 

\end{itemize}

Although the lower bound (\ref{eq:lowerbd})
does not hold {\em pathwise}, only on average,
it provides a useful guide for our pathwise experiments. 
The upper bounds (\ref{eq:upperbd}) and (\ref{eq:upperbd2})
do not apply to any numerical
experiment conducted with non-zero $B_j$, but also serve
as a useful guide to those experiments: it is clearly
undesirable to greatly exceed the error committed by simply
trusting the data. We will hence plot the lower
and upper bounds as useful comparators for
the actual error incurred in our numerical experiments below.
We note that, for the 3DVAR example from subsection
\ref{ssec:3DVar} with complete observations, 
the upper and lower bounds coincide 
in the limit $\eta \to 0$ as then $B \to 0$.
For partial observations they differ by the second term
in the upper bound.

\subsection{Experimental Setup}
\label{ssec:results}

For all the results shown we 
choose a box side of length $L=2$. 
The forcing in Eq. (\ref{eq:nse}) is taken to be $f=\nabla^{\perp}\psi$,
where $\psi=\cos(\pi k \cdot x)$ and $\nabla^{\perp}=J\nabla$ with $J$
the canonical skew-symmetric matrix, and $k=(5,5)$.  
The method used to approximate the forward model (\ref{eq:nse})
is a modification of a fourth-order Runge-Kutta method,  
ETD4RK \cite{cox2002exponential}, 
in which the Stokes semi-group is computed exactly 
by working in the incompressible Fourier
basis  $\{\psi_{k}(x)\}_{k \in {\mathbb Z}^2\backslash\{0\}}$
in Eq. (\ref{psik}),
and Duhamel's principle (variation of constants formula) is 
used to incorporate the nonlinear term.
Spatially, a Galerkin spectral method \cite{hesthaven2007spectral} 
is used, in the same basis, 
%to evaluate the nonlinear term 
and the convolutions arising from products in the nonlinear 
term are computed via FFTs.
We use a double-sized domain in each dimension, buffered with
zeros, resulting in $64^2$ grid-point FFTs, and only half the 
modes in each direction are retained when transforming back 
into spectral space in order to prevent aliasing, which is avoided
as long as fewer than 2/3 of the modes are retained.

The dimension of the attractor is determined by the 
viscosity parameter $\nu$. For the particular forcing used
there is an explicit steady state for all $\nu>0$
and for $\nu \geq 0.035$ this solution is stable (see
\cite{majda2006non}, Chapter 2 for details). 
As $\nu$ decrease the flow becomes 
increasingly complex
and the regime $\nu \leq 0.016$ corresponds to 
strongly chaotic dynamics with 
%attributes of turbulent scalings 
an upscale cascade of energy in the spectrum.  
We focus subsequent studies of the filter on a  
strongly chaotic ($\nu = 0.01$) parametric regime.
For this small viscosity parameter, we use a time-step 
of $\delta t = 0.005$.

The data is generated by computing a true signal
solving \eref{eq:nse} at the desired value of $\nu$, and
then adding Gaussian random noise to it at each observation
time. Such noise does not satisfy Assumption \ref{a:1z}, since
the supremum of the norm of the noise sequence is not finite, and so 
this setting provides a severe test beyond what is 
predicted by the theory; nonetheless, it should be noted that 
Gaussian random variables only obtain arbitrarily large values 
arbitrarily rarely. 
 
All experiments are conducted using the 3DVAR setup and
it is useful to reread the end of subsection \ref{ssec:3DVar}
in order to interpret the parameters $\alpha$ and $\eta.$ 
We consider both the choices $\alpha=\pm 1$ for 3DVAR, noting
that in the case $\alpha=-1$ the operator $B$ has norm strictly
less than one and so we expect the algorithm to be more robust
in this case (see Remark \ref{r:2} for discussion of this fact). 
For %the discrete time 
all experiments we set $\ell=\lambda_1^{-1}$
which ensures that the action of $A_0^{2\alpha}$, and hence
$B$, on the first eigenfunction is independent of the value
of $\alpha$; this is a useful normalization when comparing
computations with $\alpha=1$ and $\alpha=-1.$
%In the continuous time experiments we set $\ell=1.$

We set the observational noise to constant white
noise $\Gamma_j = \Gamma = \eps^2 I$ (i.e. $\beta=0$ in section \ref{ssec:3DVar}).  
Here $\eps = 0.04$, which gives
a standard deviation of approximately 10\% of the maximum 
standard deviation of the strongly chaotic dynamics.  
Since we are computing in a truncated finite-dimensional basis 
the eigenvalues are summable; the situation 
can be considered as an approximation of an operator
whose eigenvalues decay rapidly outside the basis in which
we compute. 

To be more precise regarding the algorithm, let $U$ denote the 
finite-dimensional spectral representation, which is a complex-valued
vector of dimension $32^2$, including redundancy arising from reality 
and zero mass constraints.  The computation of $\cB$ in 
Eq. (\ref{eq:nse}) requires padding this with zeros to a $64^2$ 
vector, computing inverse FFTs on the discretization of the 
6 fields $u_i, u_{i,j}$ for $i,j \in \{1,2\}$, performing products, 
computing FFTs on the 2 resulting (discrete) spatial fields, and
finally discarding the now-populated padding modes.
Denoting the discrete map from time $t$ to time $t+s$ 
by $\Phi_s(M)$, the experiments proceed precisely as follows:
\begin{itemize}
\item Evolve $U_t=\Phi_t(U_0)$ until statistical equilibrium, as judged
by observing the energy fluctuations $E[U_t(t)] = ||U_t||_2^2$.  
Set $U_0=U_t$ so that the initial condition is on the attractor.
\item Compute the observations $Y_j = \Phi_{j h}(U_0) + \cN(0,\Gamma)$.
\item Draw $\hU_0 \sim \cN(0,\kappa A^{2\alpha})$, where $\kappa \gg
  1$ [this is essentially arbitrary, as long as the initial condition
  is something sensible and such that $\| \hU_0-U_0\| = O(1)$].
\item Compute $B$, and $I-B$ as given in Equations (\ref{be}, \ref{ibe})
\item For $j=1, \dots, J$ : Compute $\hU_j = B \Phi_h(\hU_{j-1}) + (I-B) Y_j.$

\end{itemize}

\subsection{Complete Observations}
\label{ssec:complete}

\begin{figure*}
\includegraphics[width=.45\textwidth]{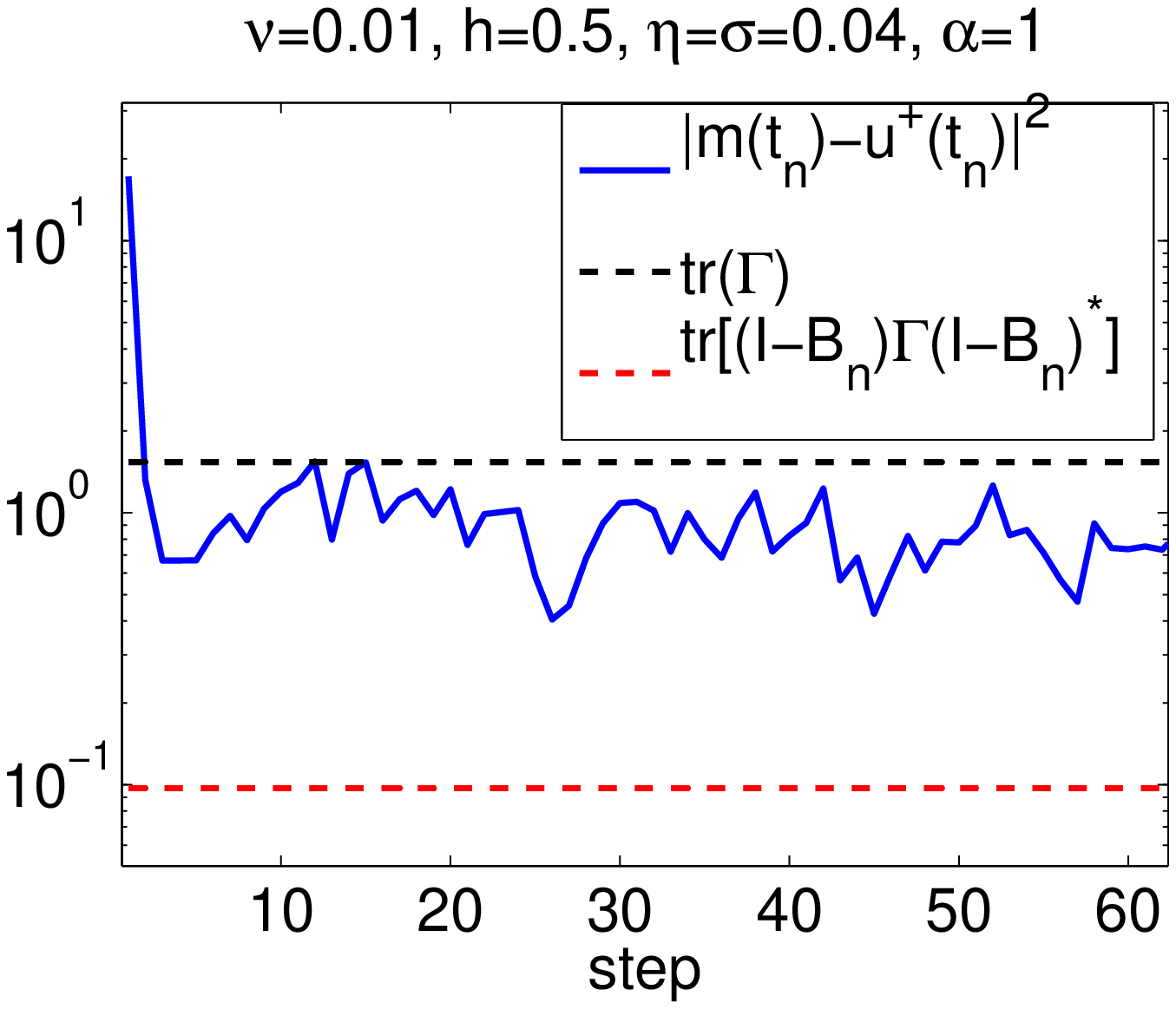}
\includegraphics[width=.45\textwidth]{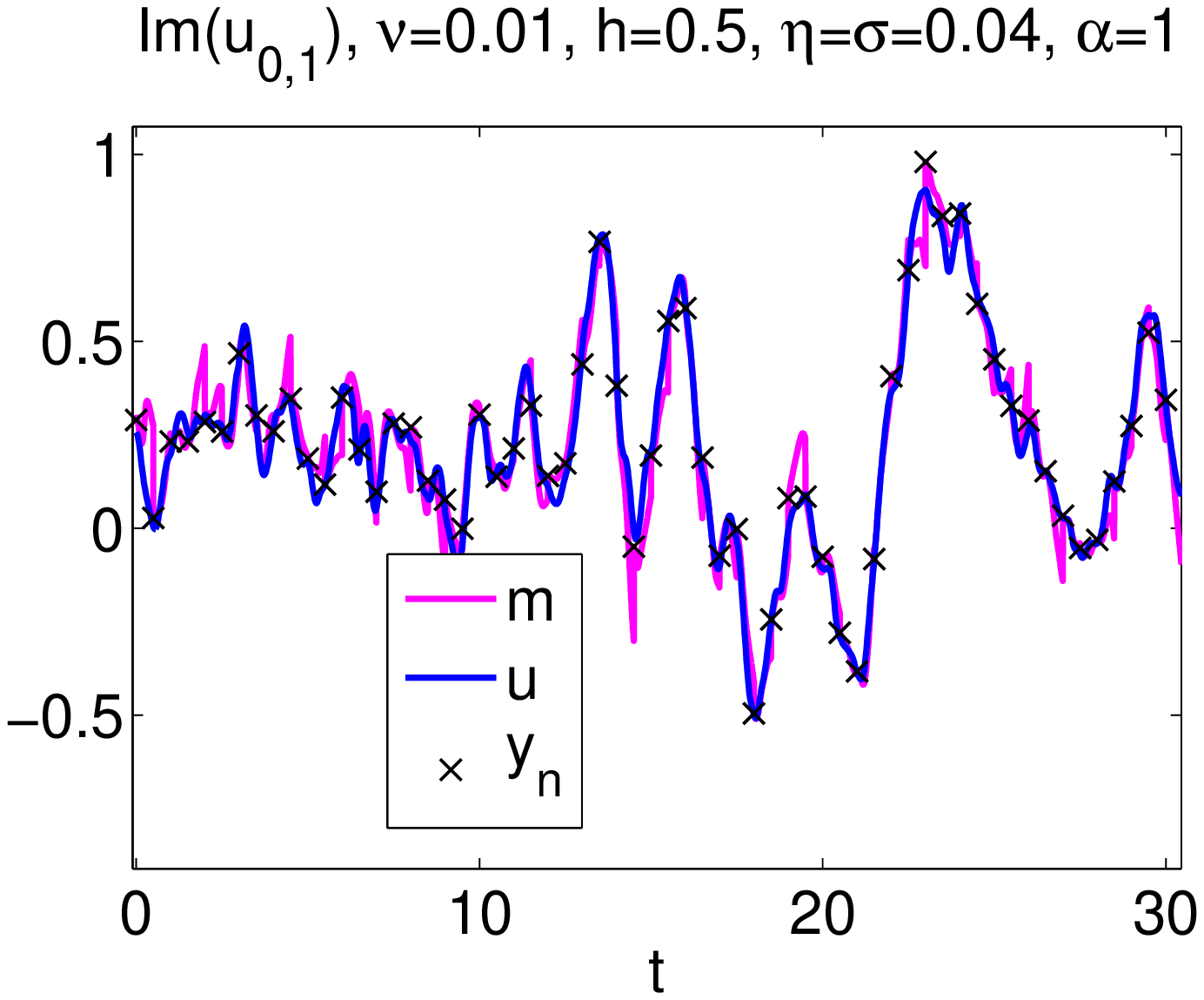}
\includegraphics[width=.45\textwidth]{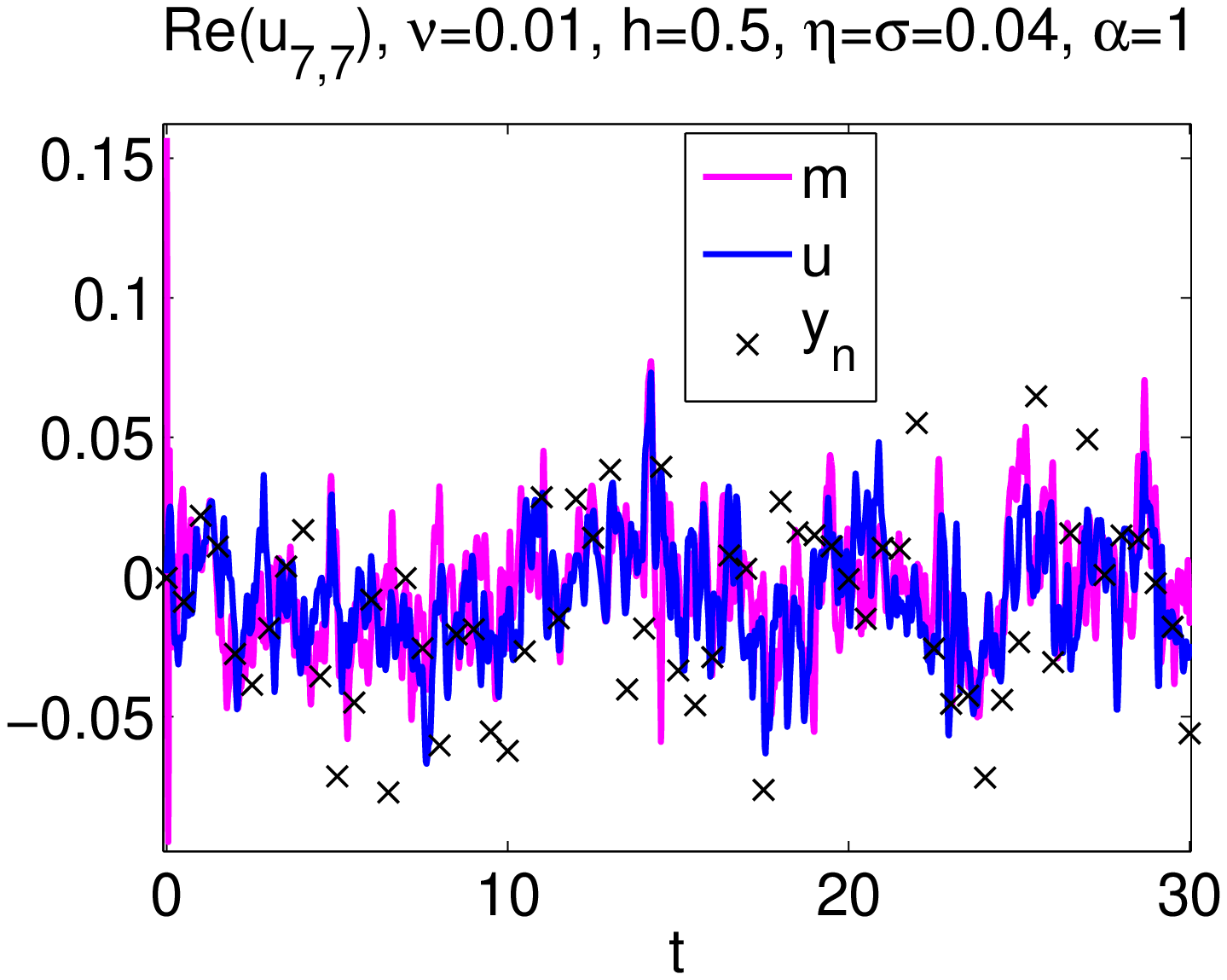}
\includegraphics[width=.45\textwidth]{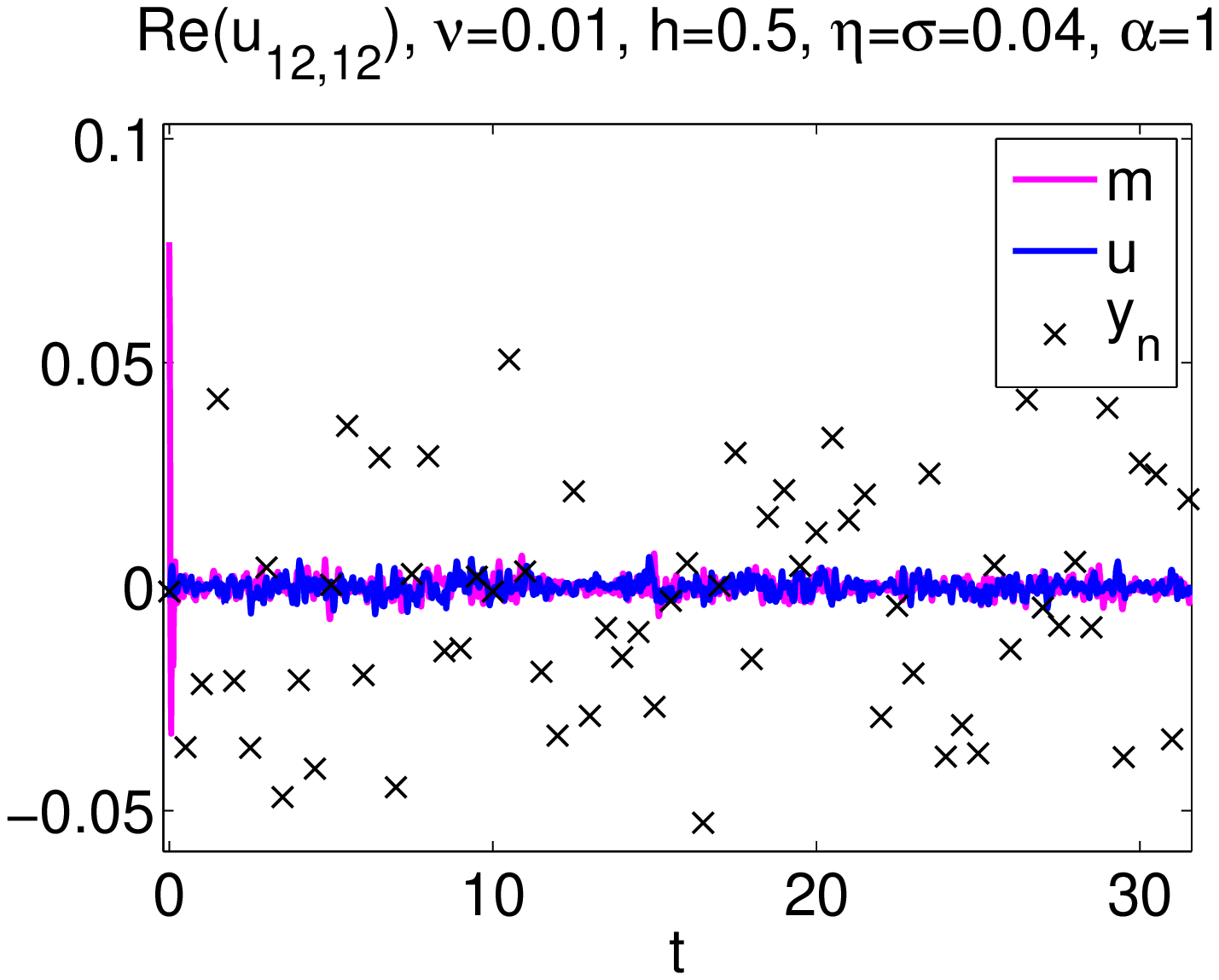}
\caption{Example of a stable trajectory for 3DVAR with $\nu=0.01,
  h=0.5, \eta=\sigma=0.04, \alpha=1$.  
The top left plot shows the norm-squared error between the
estimated mean, $m(t_n)=\hat{m}_n$, and the signal, $u(t_n)$, 
in comparison to the preferred upper
bound (i.e. the total observation error ${\rm tr} (\Gamma) = \Xi$) 
and the lower bound 
${\rm tr} [(I-B_n) \Gamma (I-B_n)]$.  The other three plots show the 
estimator, $m(t)$, together with the signal, $u(t)$, and the
observations, $y_n$ for a few individual modes.}
\label{a1.1}
\end{figure*}

\begin{figure*}
\includegraphics[width=.45\textwidth]{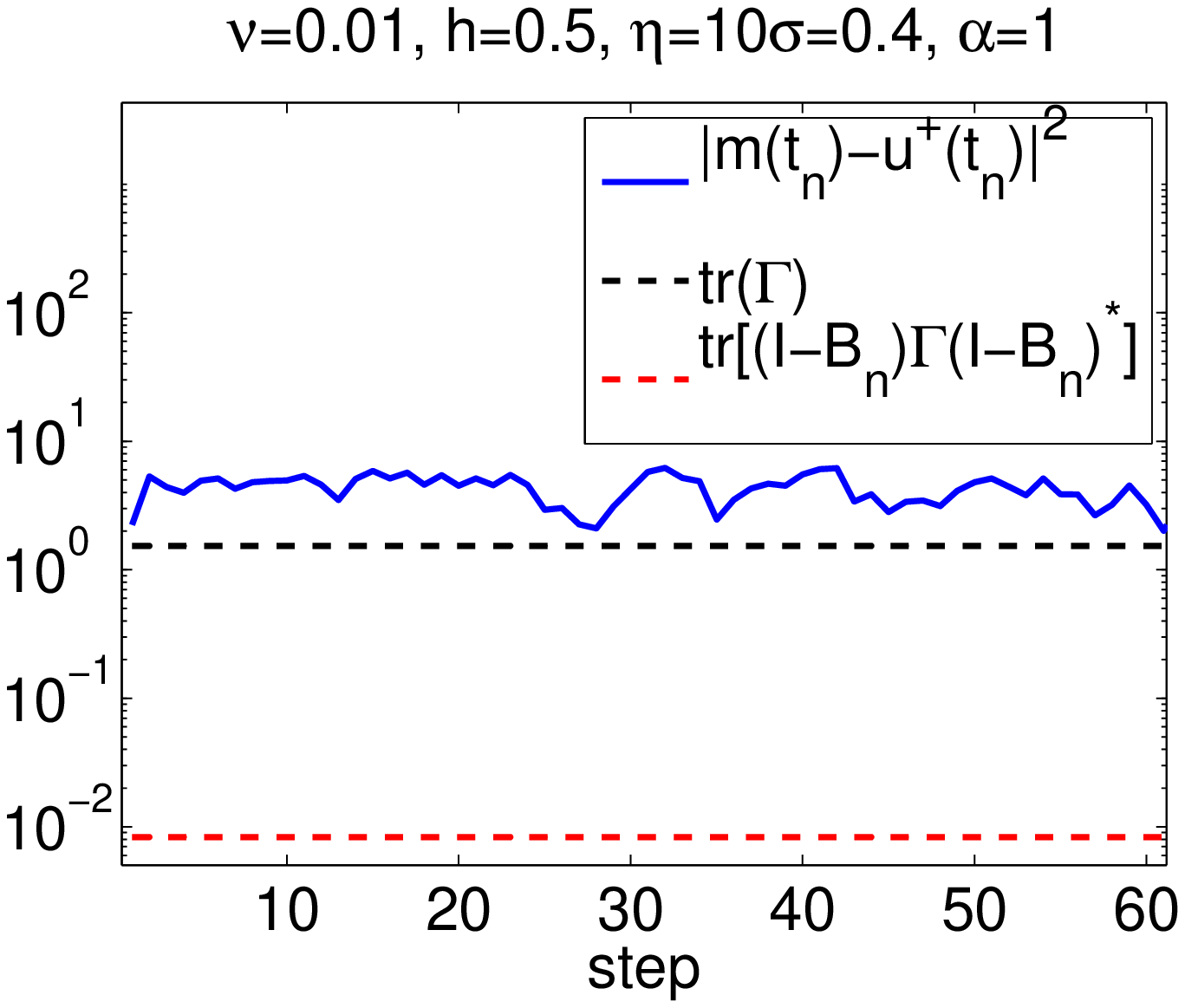}
\includegraphics[width=.45\textwidth]{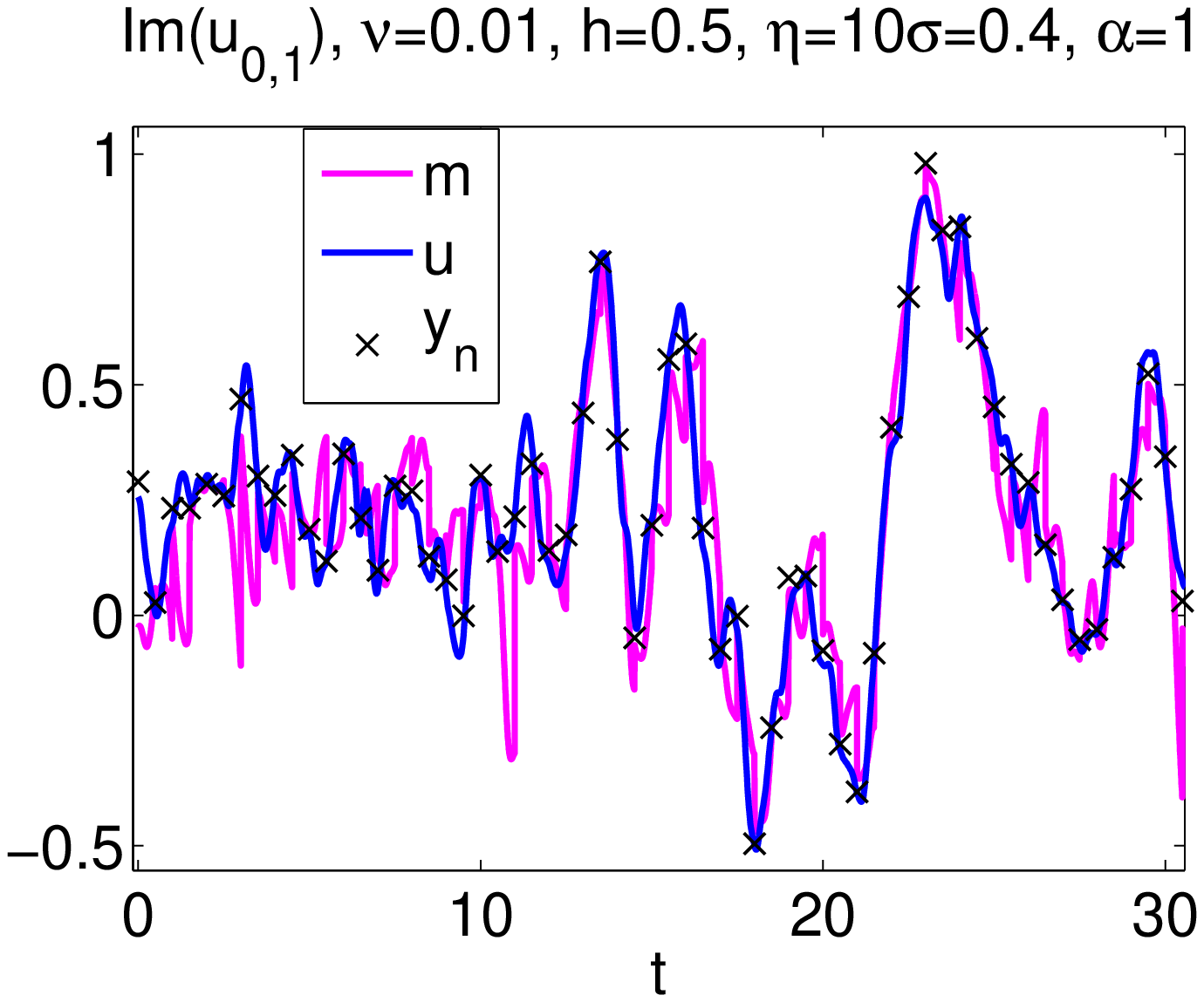}
\includegraphics[width=.45\textwidth]{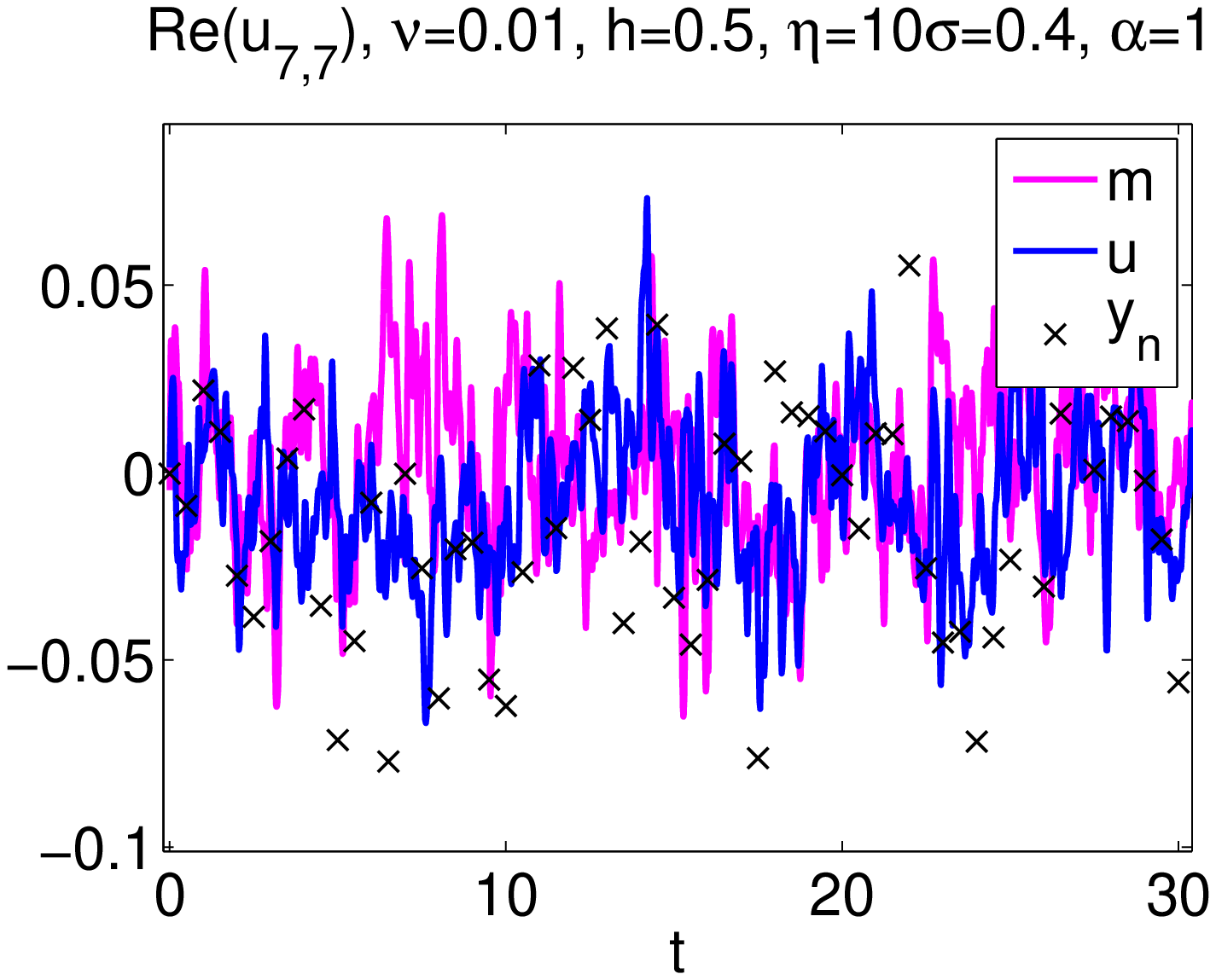}
\includegraphics[width=.45\textwidth]{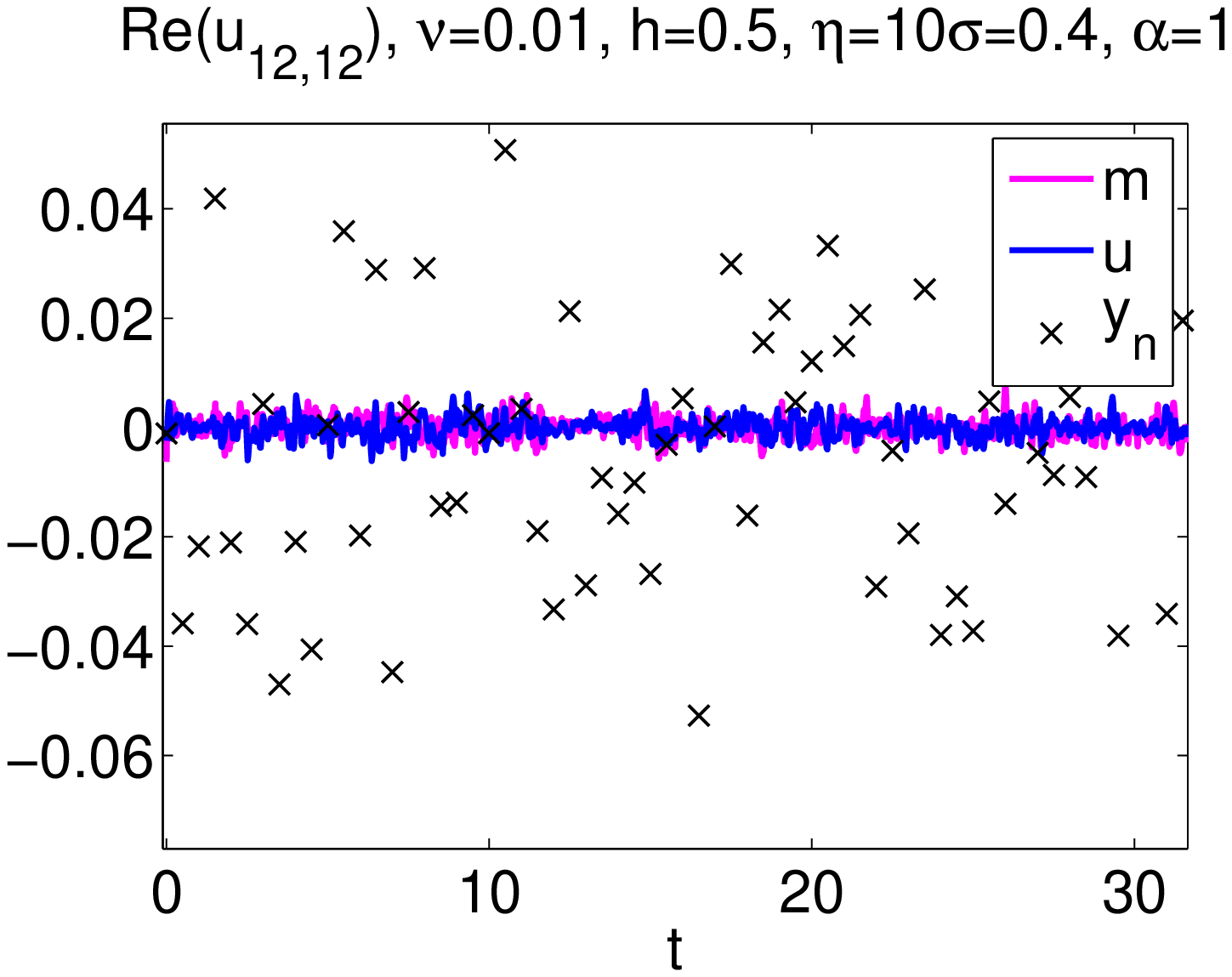}
\caption{Example of a destabilized trajectory for 3DVAR with 
the same parameters as in Fig. \ref{a1.1} except the larger
value of $\eta=10\sigma=0.4$.  Panels are the same.}
\label{a1.10}
\end{figure*}

\begin{figure*}
\includegraphics[width=.45\textwidth]{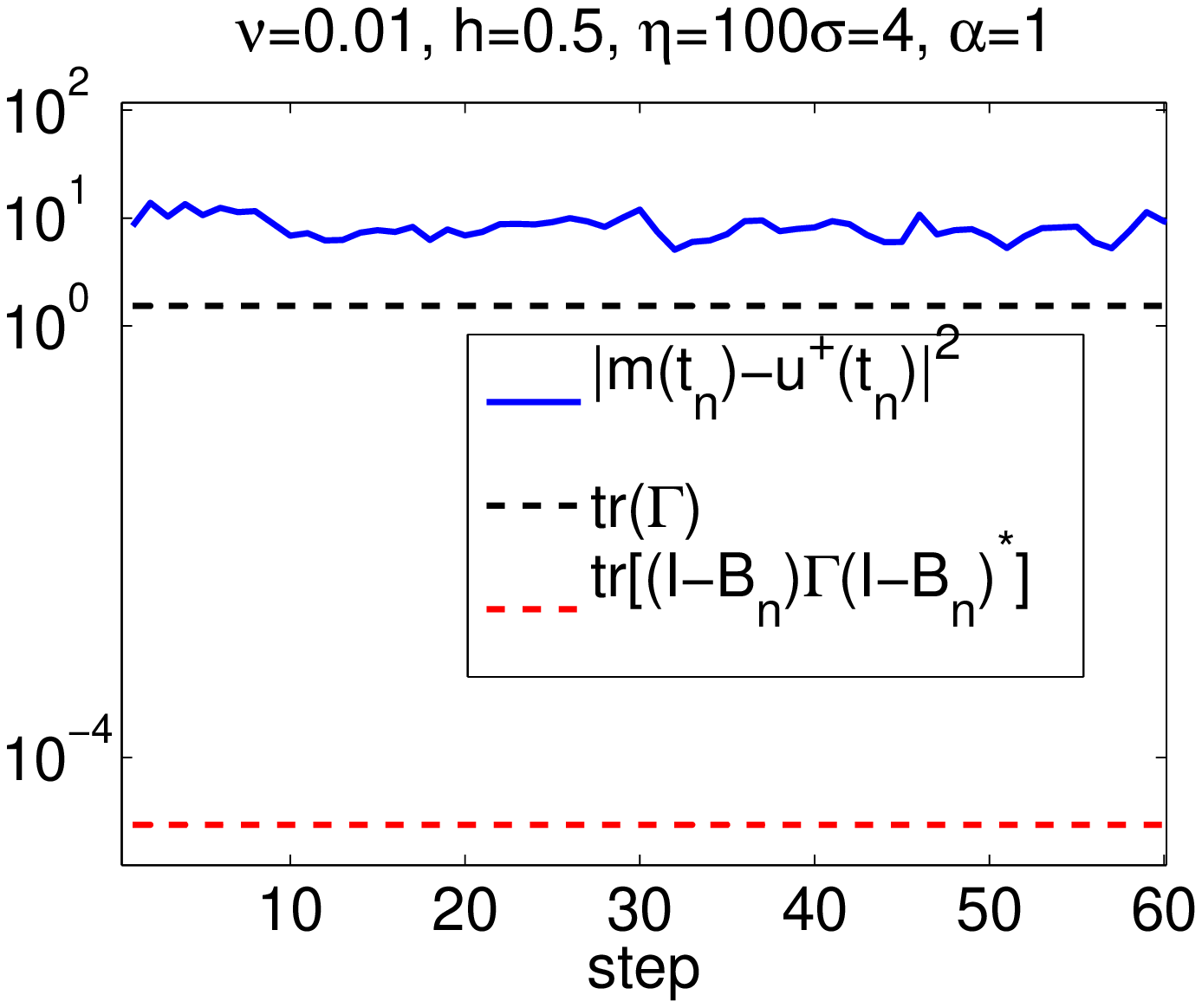}
\includegraphics[width=.45\textwidth]{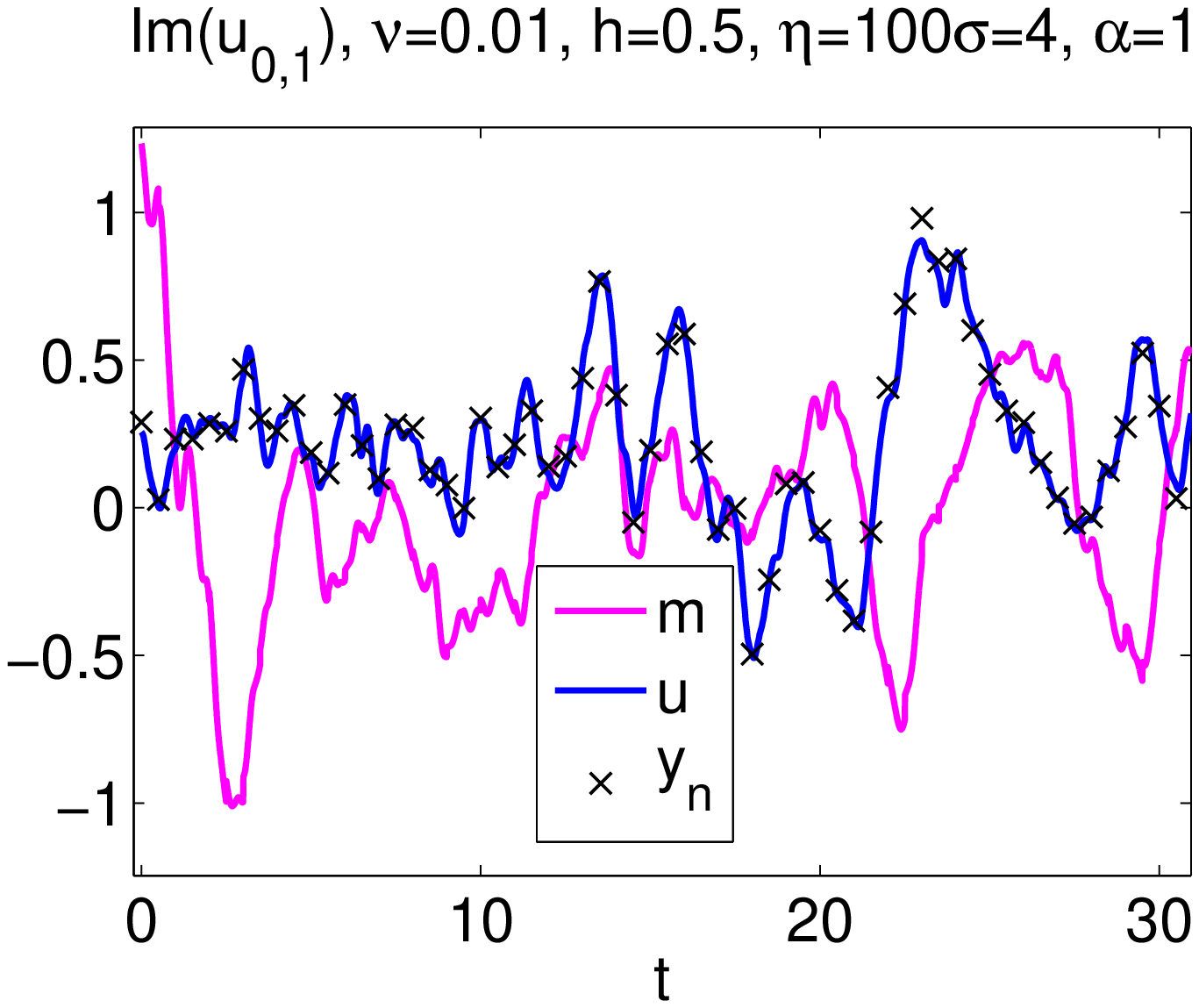}
\includegraphics[width=.45\textwidth]{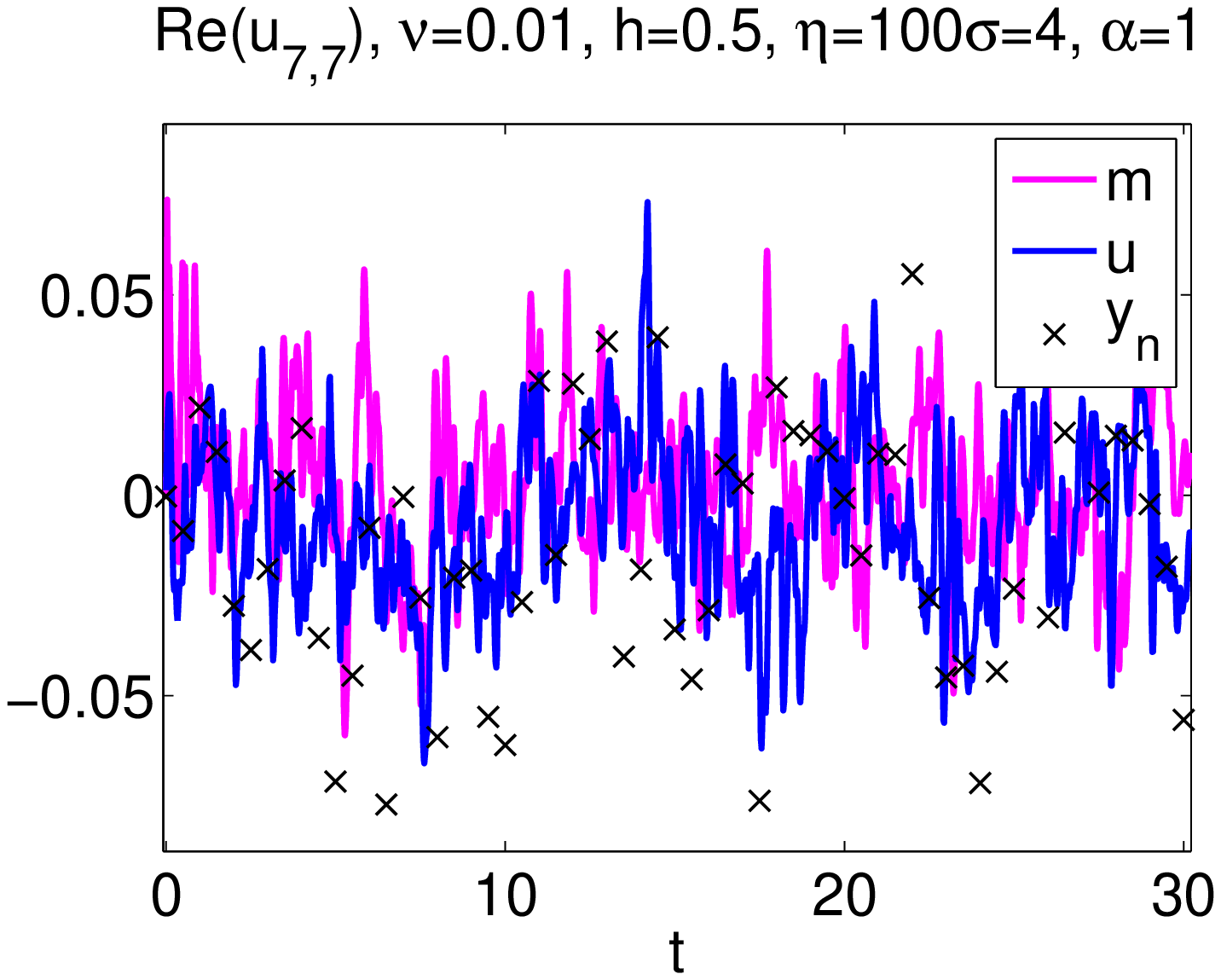}
\includegraphics[width=.45\textwidth]{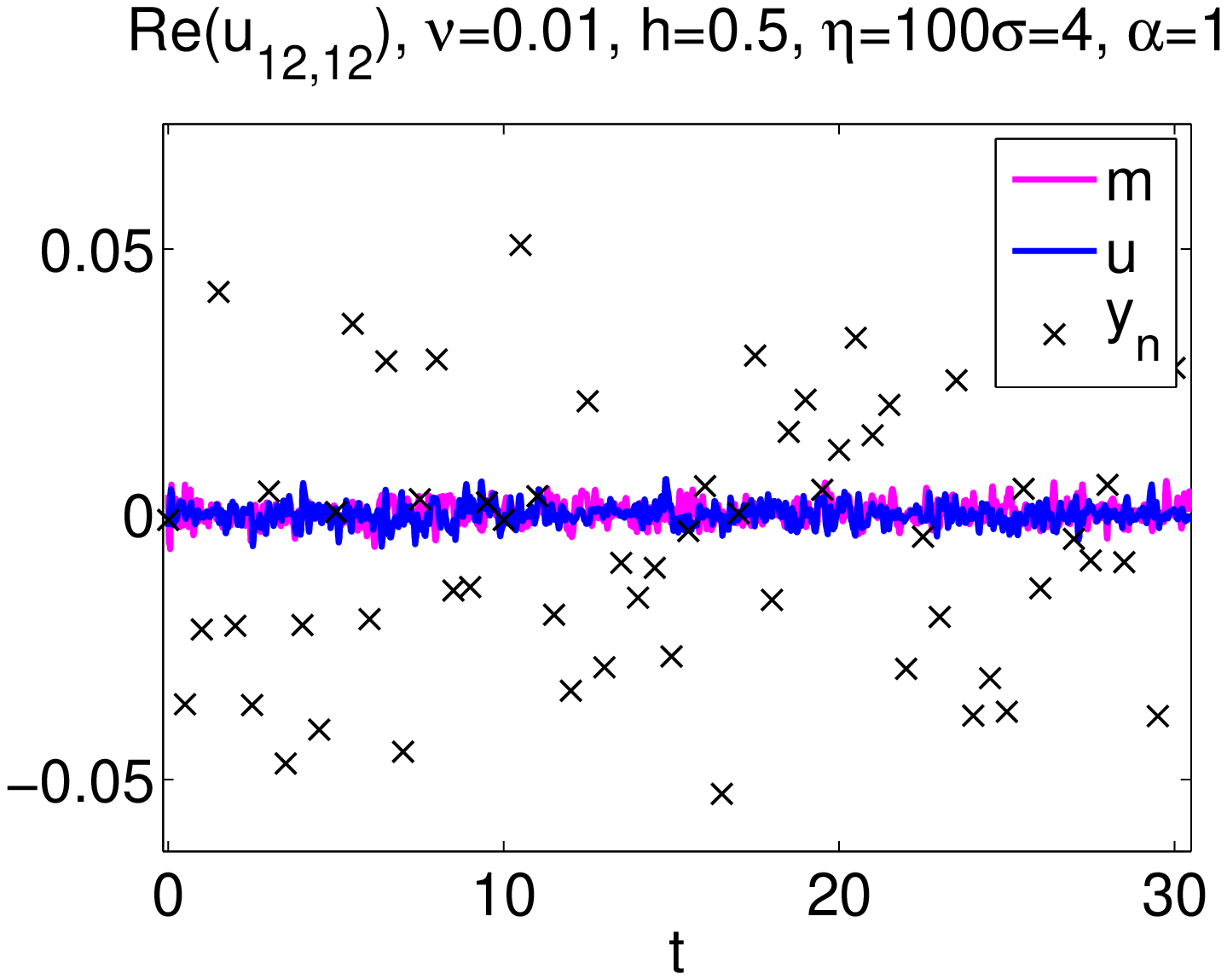}
\caption{Example of a destabilized trajectory for 3DVAR with 
the same parameters as in Fig. \ref{a1.1} except the larger
value of $\eta=100\sigma=4$.  Panels are the same.}
\label{a1.100}
\end{figure*}

We start by considering discrete and complete observations and
illustrate Theorem \ref{t:mz}, and in particular the role of 
the parameter $\eta.$
The experiments presented employ a
large observation increment of $h = 0.5 = 100 \delta t$.
For $\alpha=1$ we find that 
when $\eta=\sigma$ (Fig. \ref{a1.1}) the estimator 
stabilizes from an initial ${\mathcal O}(1)$ error
and then remains stable. The upper and lower bounds are
satisfied (the upper bound after an initial rapid transient), 
and even the high modes, which are slaved to the low modes, 
synchronize to the true signal. 
For $\eta=10\sigma$ (Fig. \ref{a1.10})
the estimator fails to satisfy the upper bound, 
but remains stable over a long time horizon; there
is now significant error in the $k=(7,7)$ mode,
in contrast to the situation with smaller $\eta$ shown 
in Fig. \ref{a1.1}.  
Finally, when $\eta=100\sigma$ (Fig. \ref{a1.100}), the estimator really
diverges from the signal, although still remains
bounded. 
 
\begin{figure*}
\includegraphics[width=.45\textwidth]{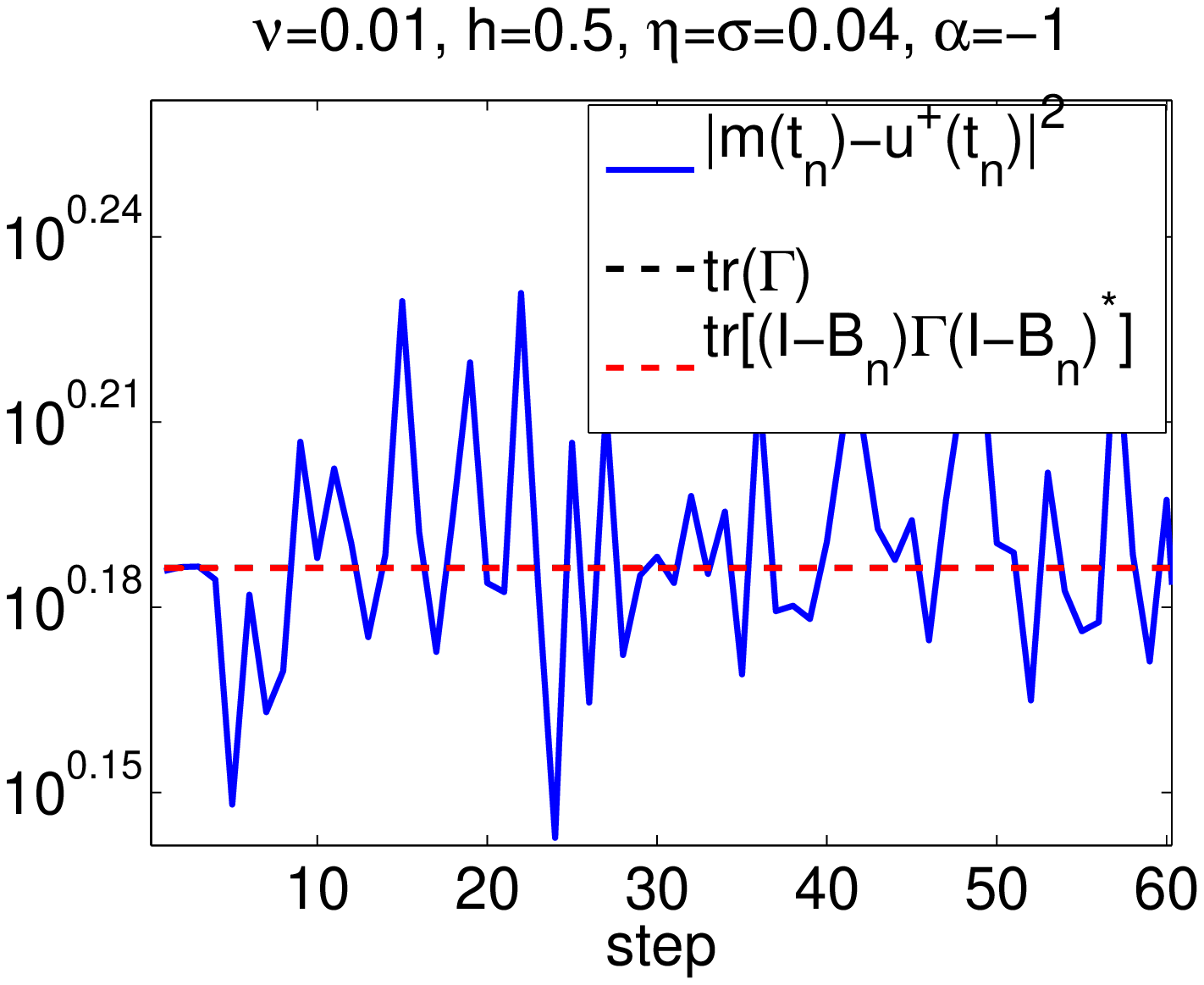}
\includegraphics[width=.45\textwidth]{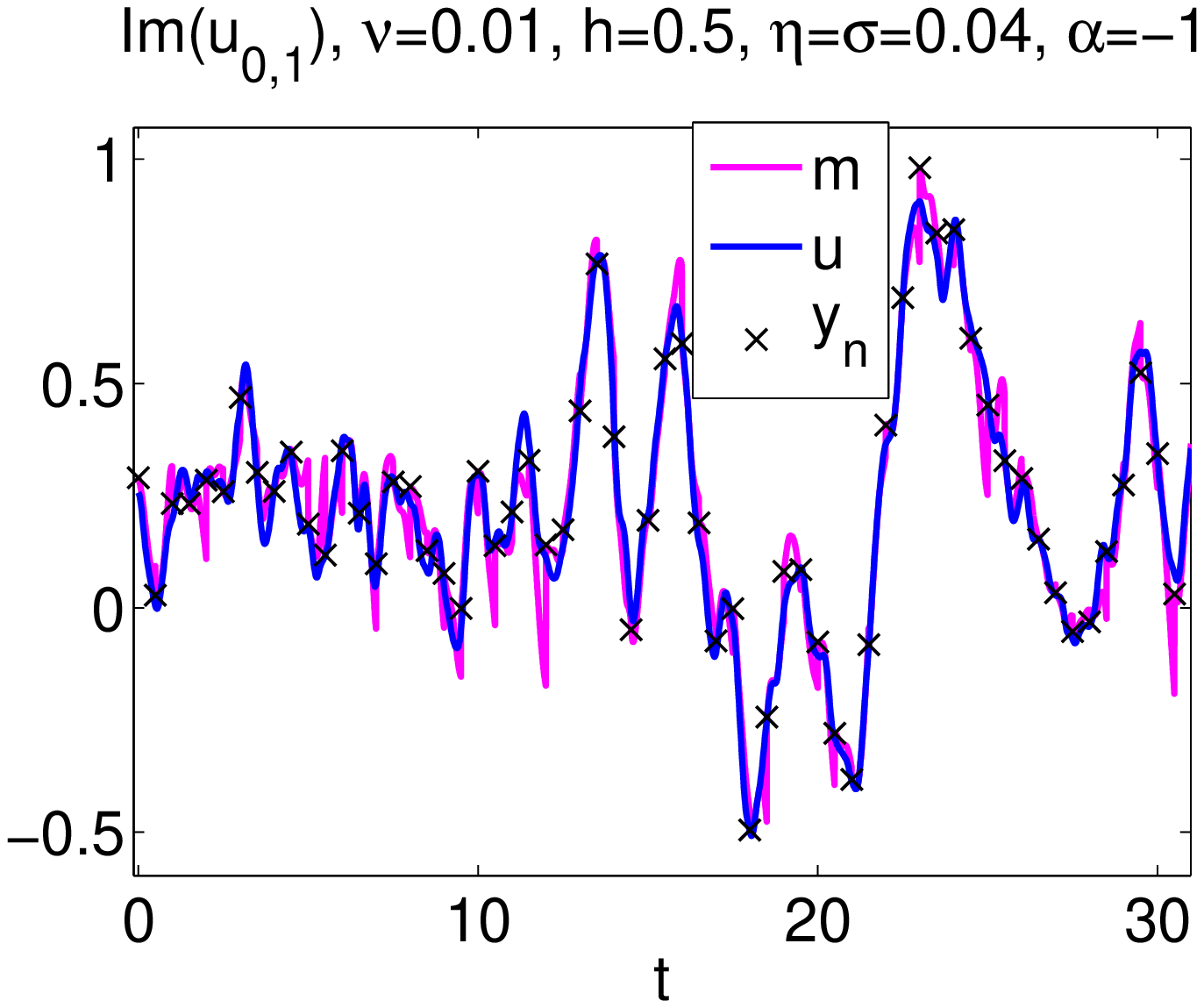}
\includegraphics[width=.45\textwidth]{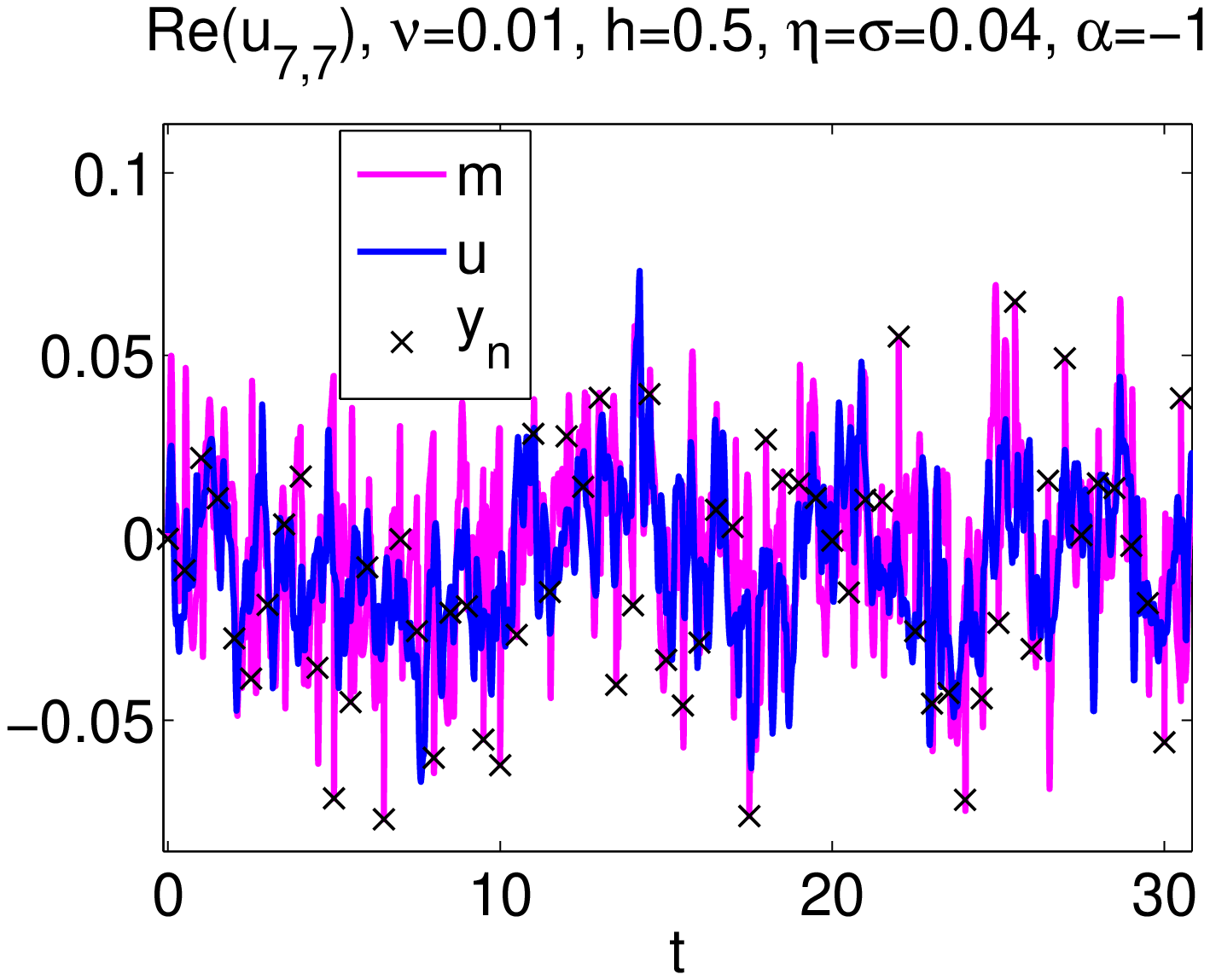}
\includegraphics[width=.45\textwidth]{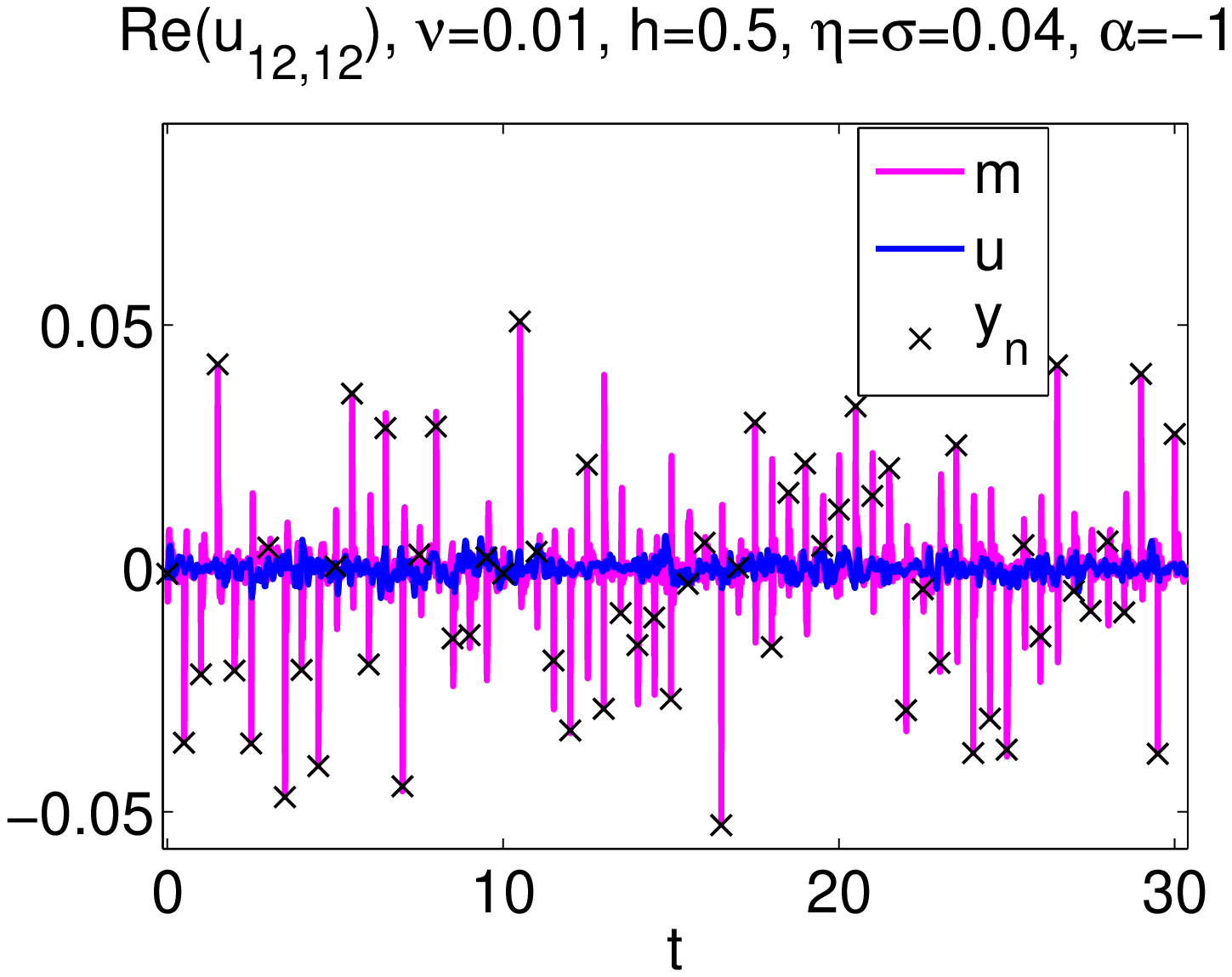}
\caption{Example of a stable trajectory for 3DVAR with 
the same parameters as in Fig. \ref{a1.1} except with
$\alpha=-1$.  Panels are the same.}
\label{am1.1}
\end{figure*}

\begin{figure*}
\includegraphics[width=.45\textwidth]{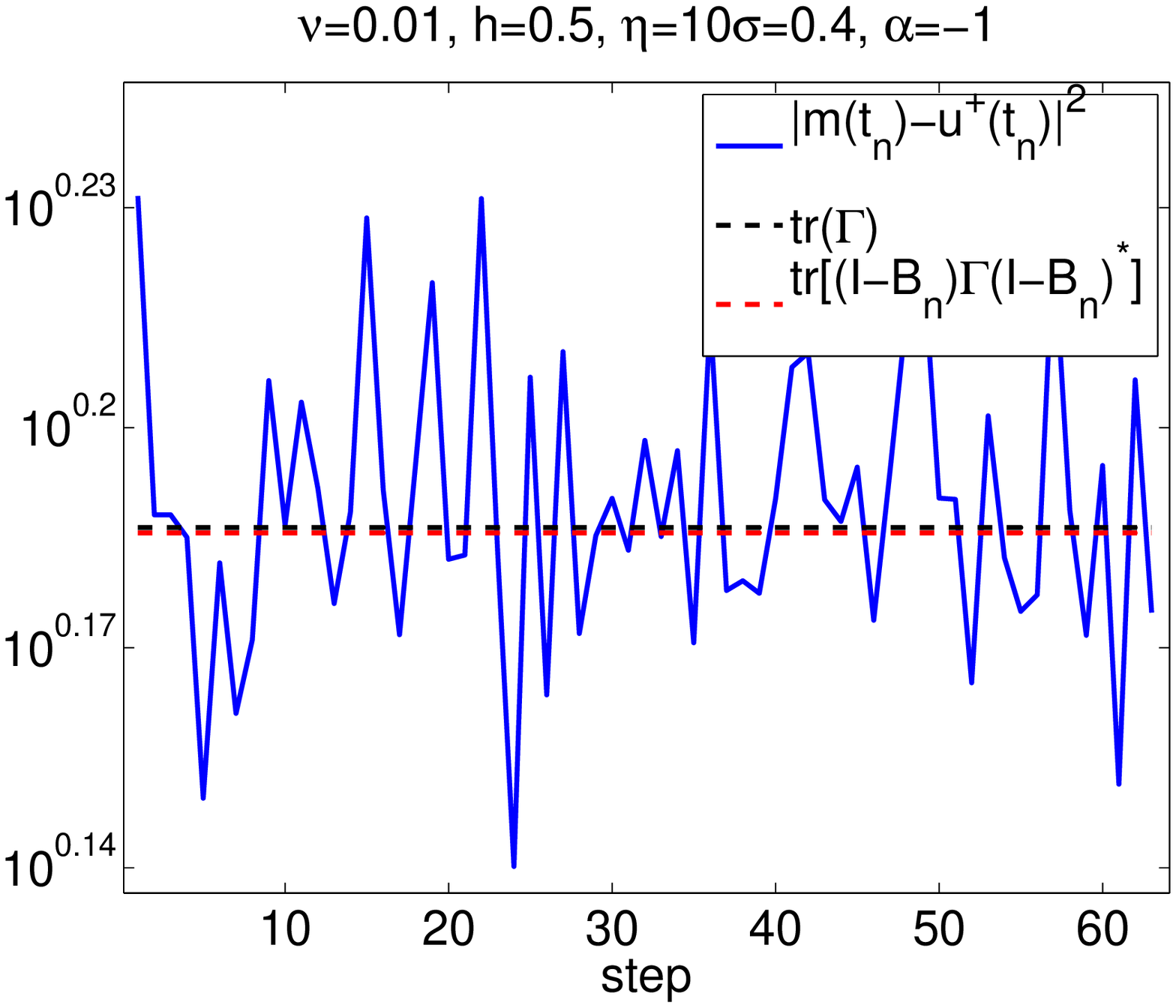}
\includegraphics[width=.45\textwidth]{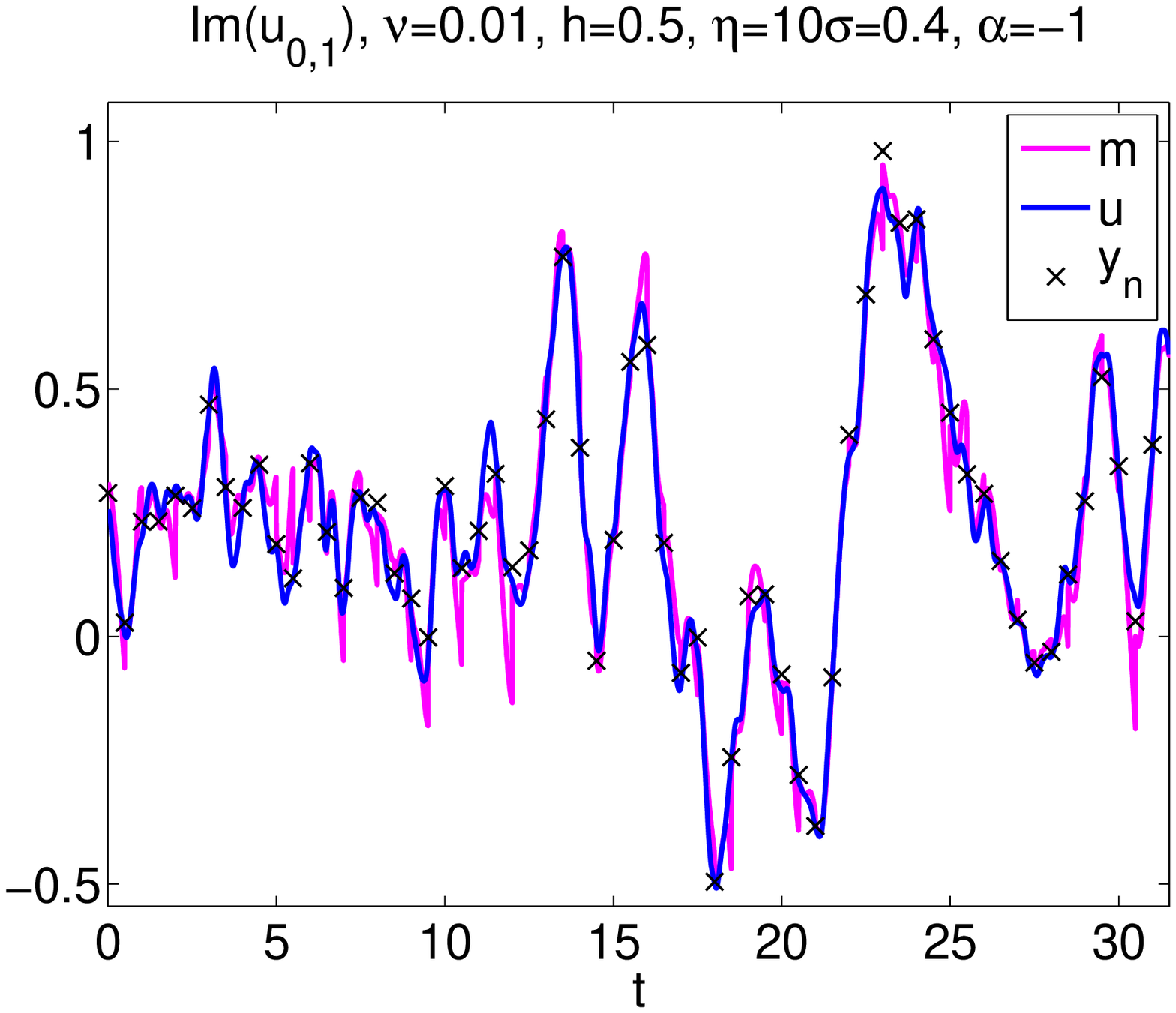}
\includegraphics[width=.45\textwidth]{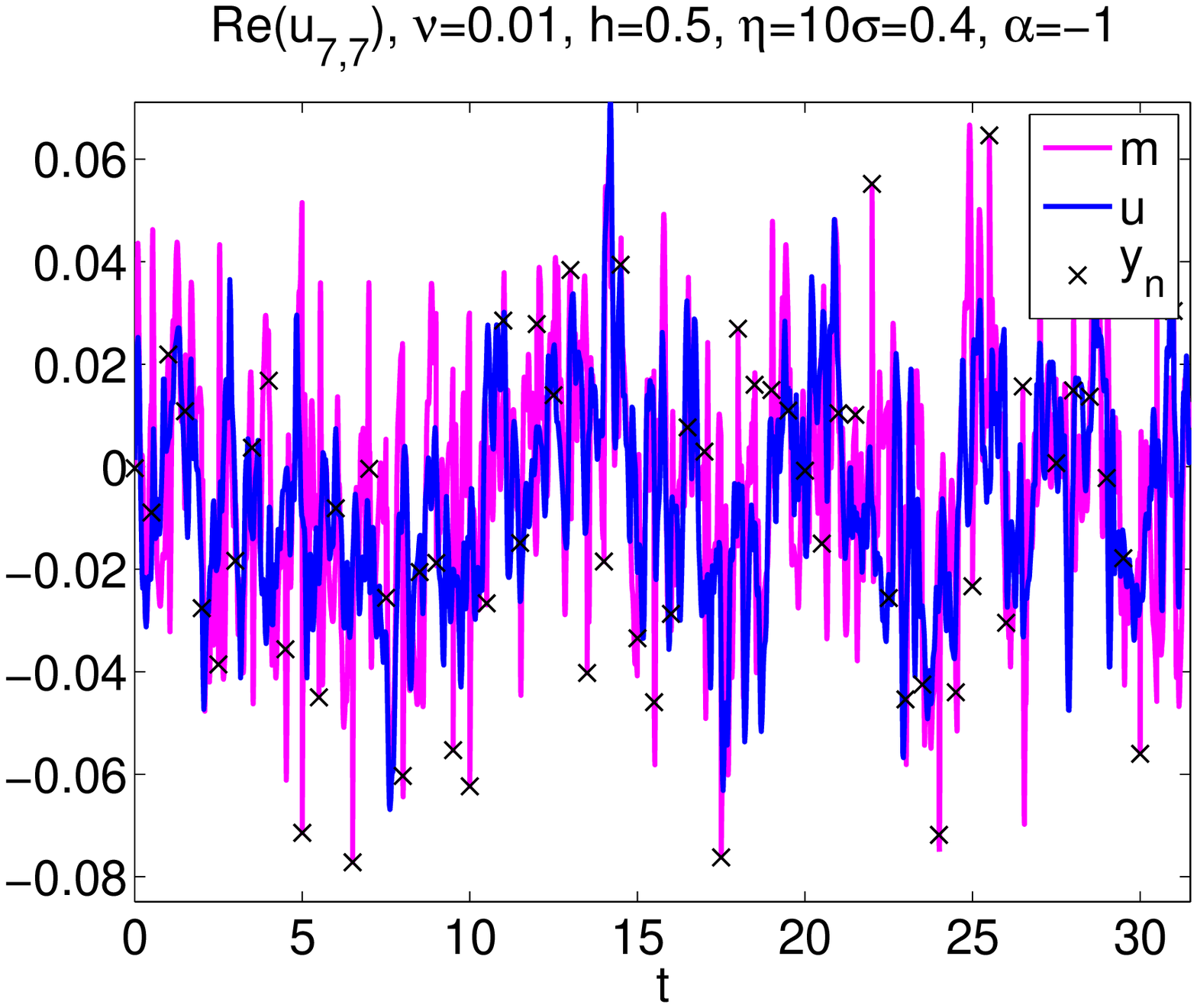}
\includegraphics[width=.45\textwidth]{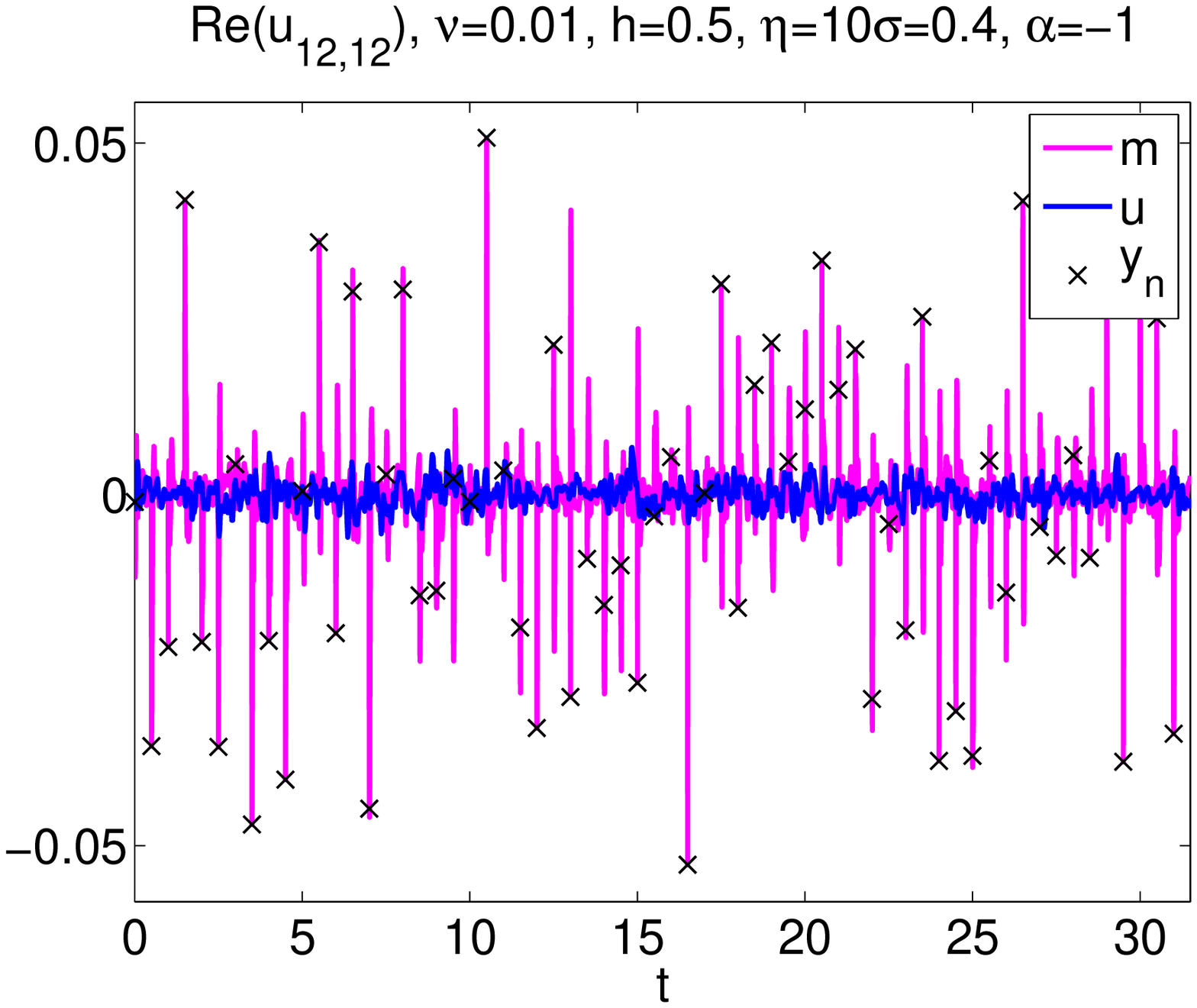}
\caption{Example of a stable trajectory for 3DVAR with 
the same parameters as in Fig. \ref{a1.10} except with
value of $\alpha=-1$.  Panels are the same.}
\label{am1.10}
\end{figure*}

\begin{figure*}
\includegraphics[width=.45\textwidth]{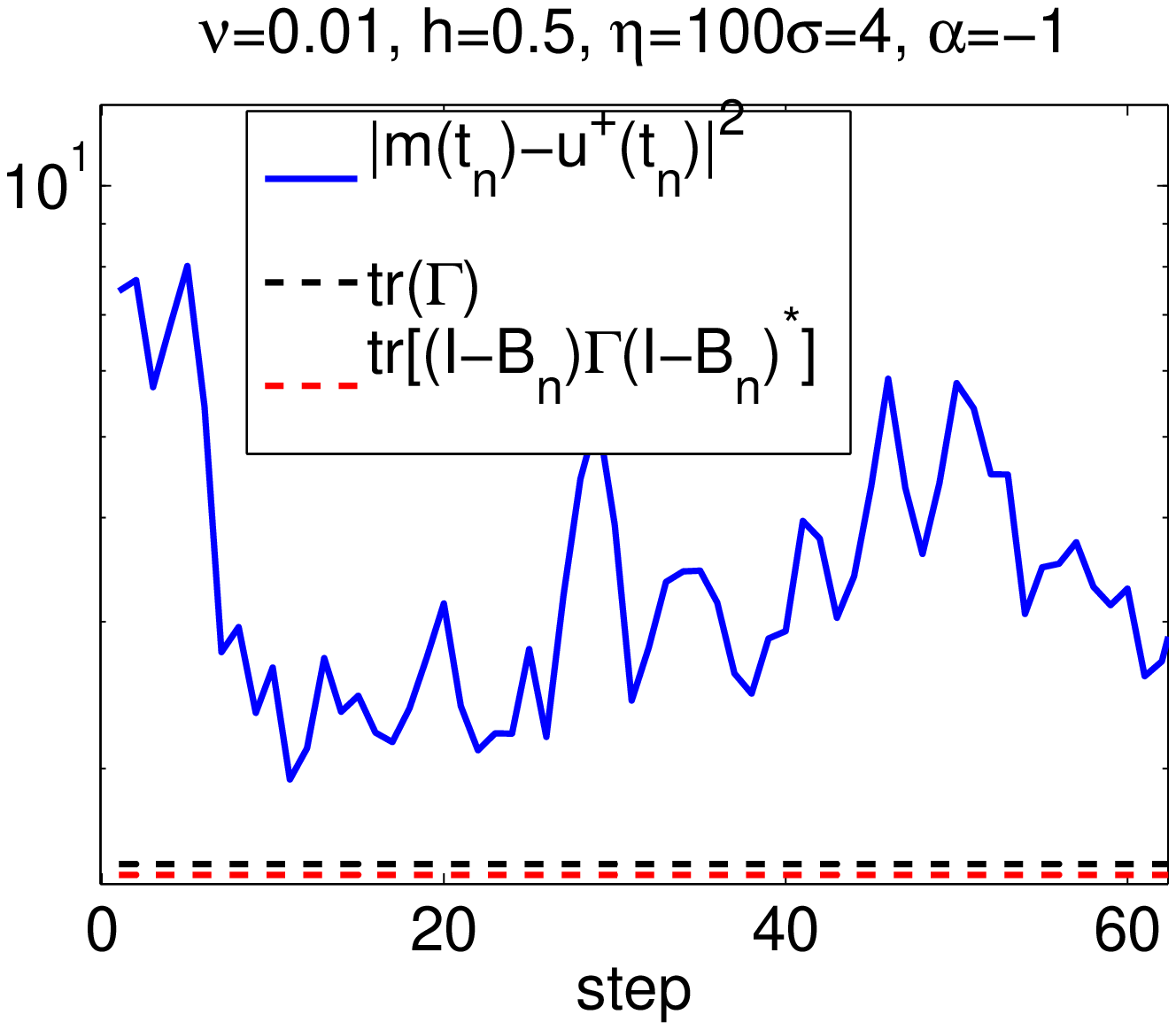}
\includegraphics[width=.45\textwidth]{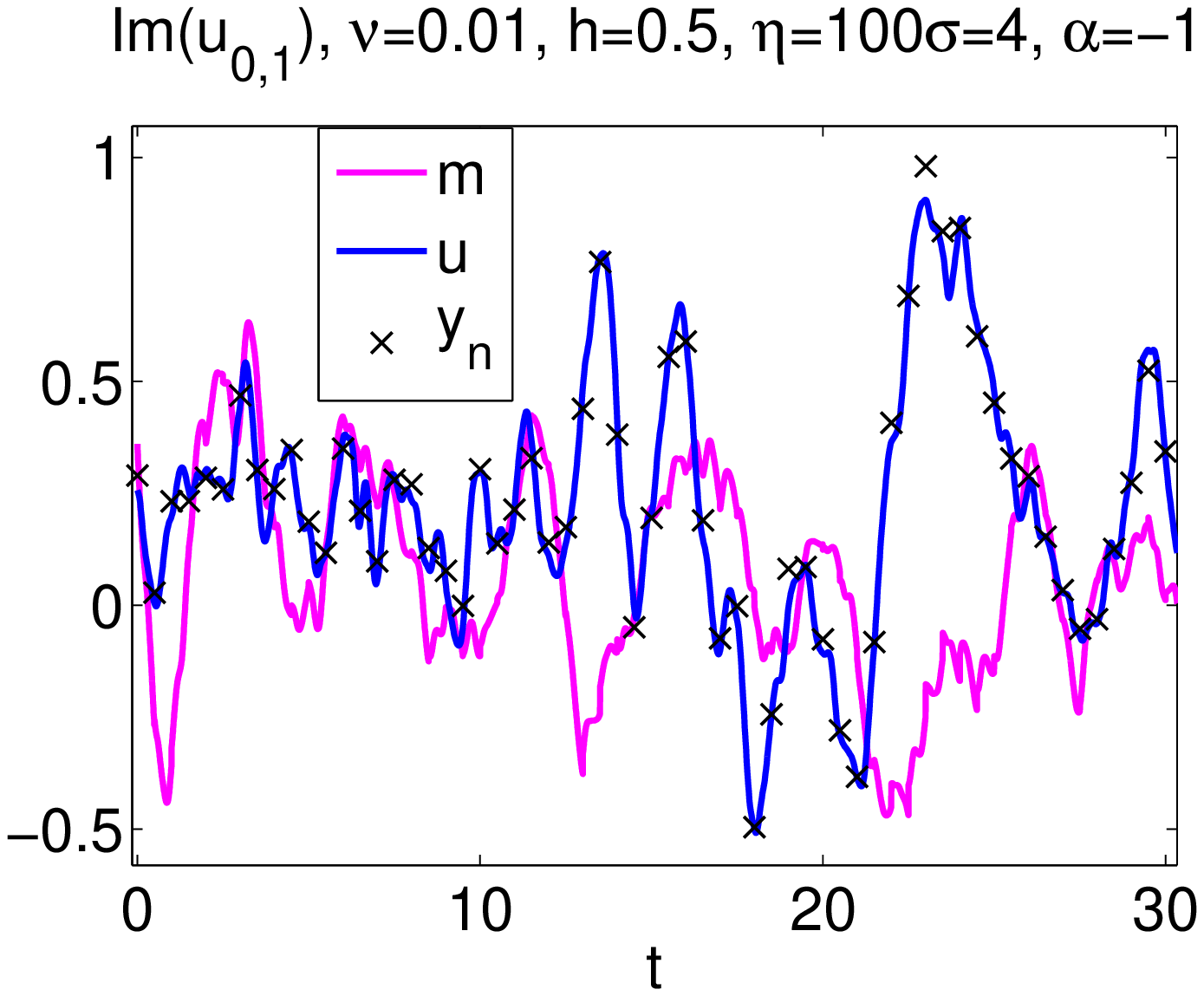}
\includegraphics[width=.45\textwidth]{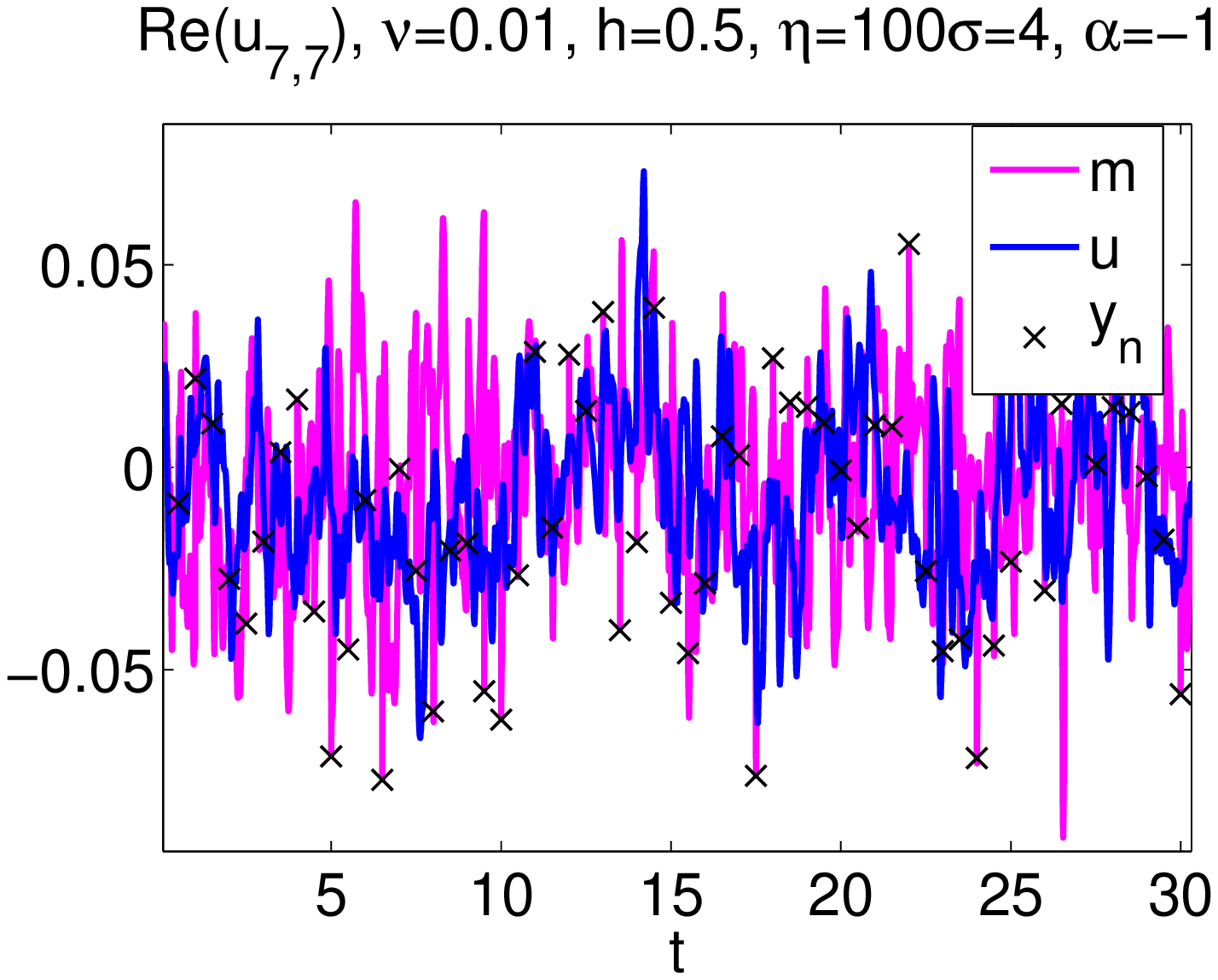}
\includegraphics[width=.45\textwidth]{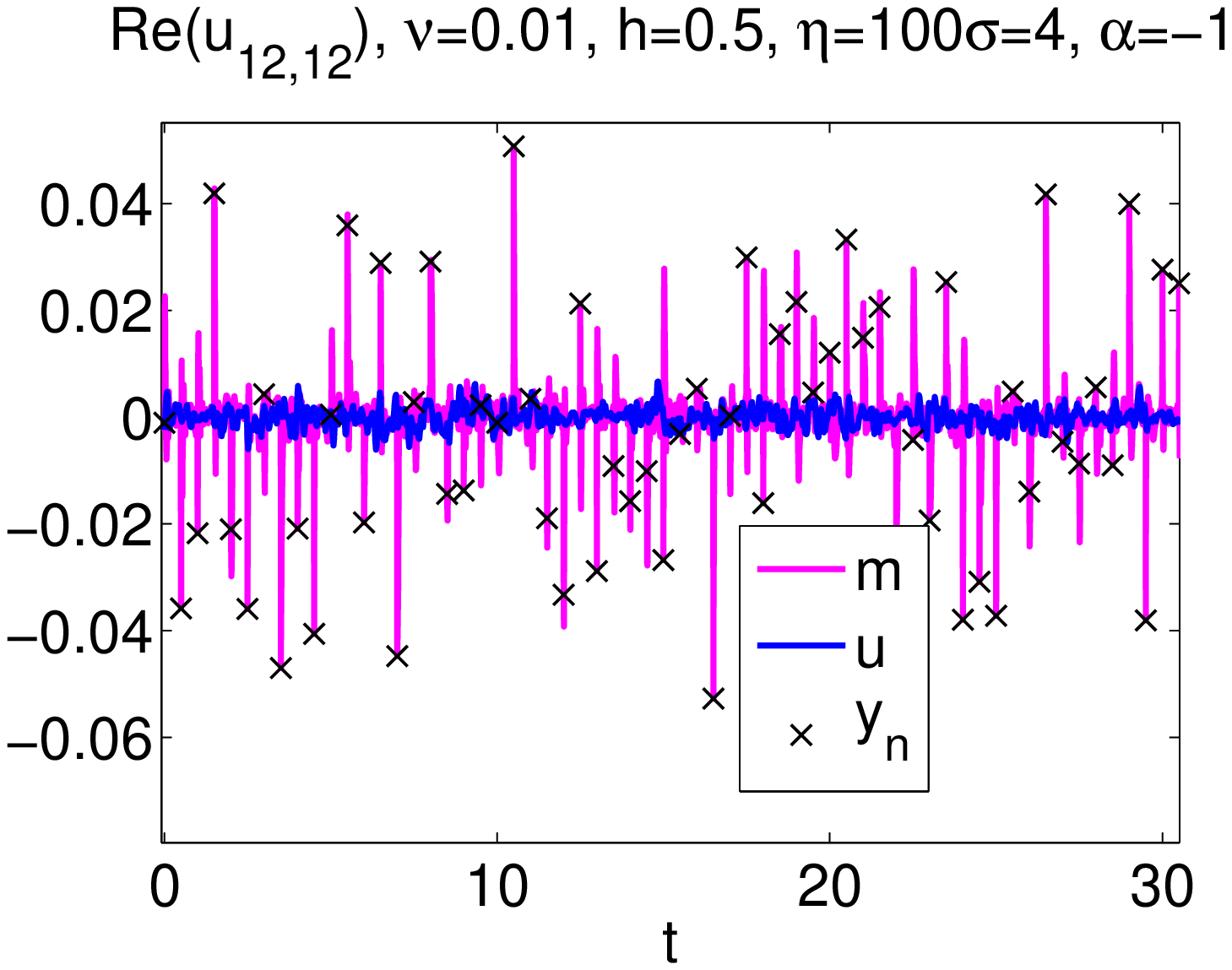}
\caption{Example of a destabilized trajectory for 3DVAR with 
the same parameters as in Fig. \ref{a1.100} except with
value of $\alpha=-1$.  Panels are the same.}
\label{am1.100}
\end{figure*}

When $\alpha=-1$ the lower and upper bounds are 
almost indistinguishable and, for all values of
$\eta$ examined, the error either exceeds or fluctuates
around the upper bound; see Figures \ref{am1.1}, \ref{am1.10}
and \ref{am1.100} where $\eta=\sigma, 10\sigma$ and $100\sigma$
respectively. It is not until $\eta=100\sigma$ (Fig. 
\ref{am1.100}) that the estimator really loses the signal.  
Notice also that the high modes of the estimator 
always follow the noisy observations and this could be 
undesirable.  For both $\eta=100\sigma$ and $10\sigma$, the $\alpha=-1$ estimator 
performs better than the one for $\alpha=1$ in terms of overall error, 
illustrating the robustness alluded to in Remark \ref{r:2}
since for $\alpha<0$ we have $\|B\|<1.$  
However, an appropriately 
tuned $\alpha=1$ filter has the potential to perform 
remarkably well, both in terms of
overall error and individual error of all modes
(see Fig. \ref{a1.1}, in contrast to 
Fig. \ref{am1.1}).  In particular,
this filter has an expected error substantially smaller than the 
upper bound, which does not happen 
for the case of $\alpha=-1$ when complete observations
are assimilated.

\subsection{Partial Observations}
\label{ssec:partial}

\begin{figure*}
\includegraphics[width=.45\textwidth]{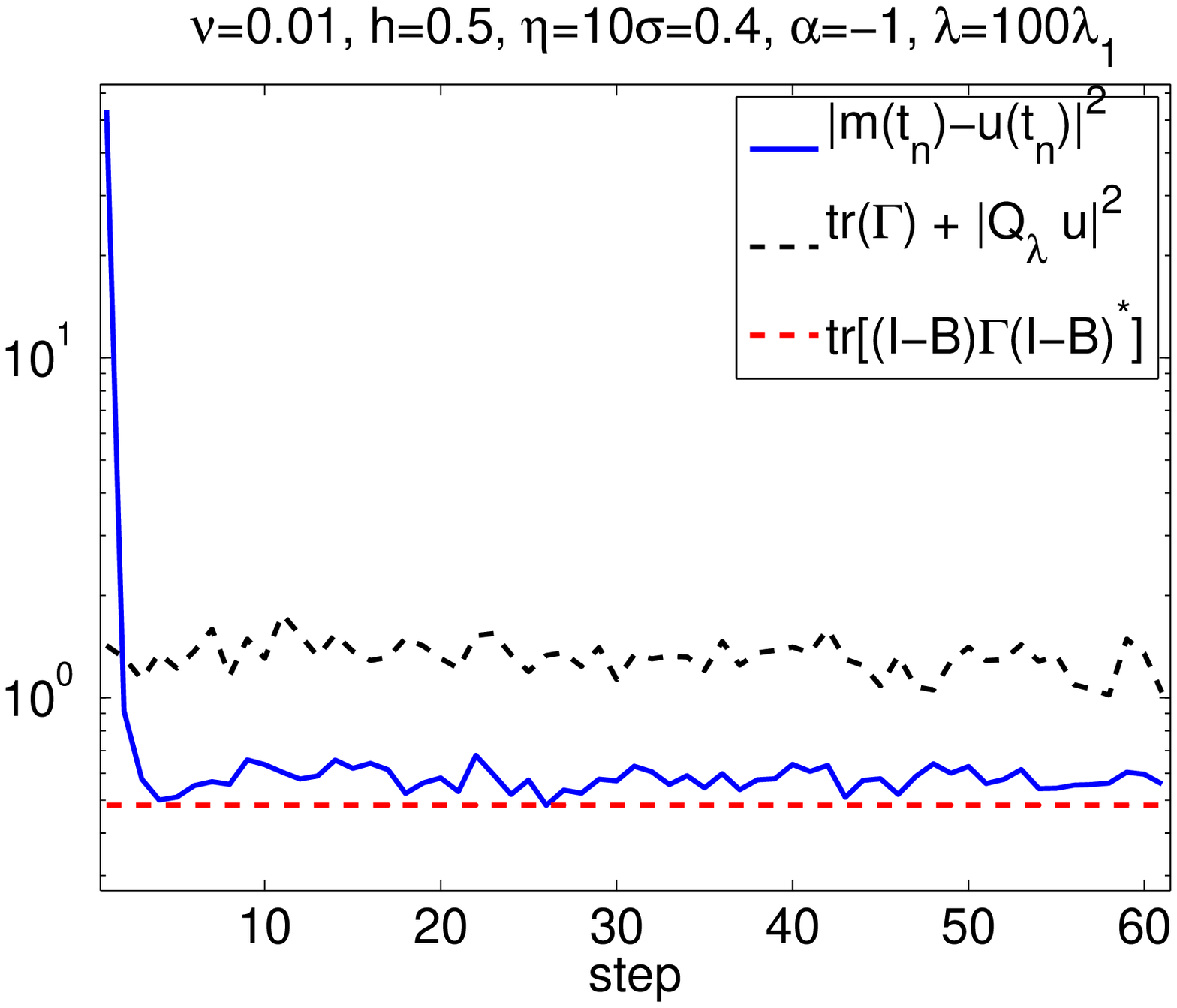}
\includegraphics[width=.45\textwidth]{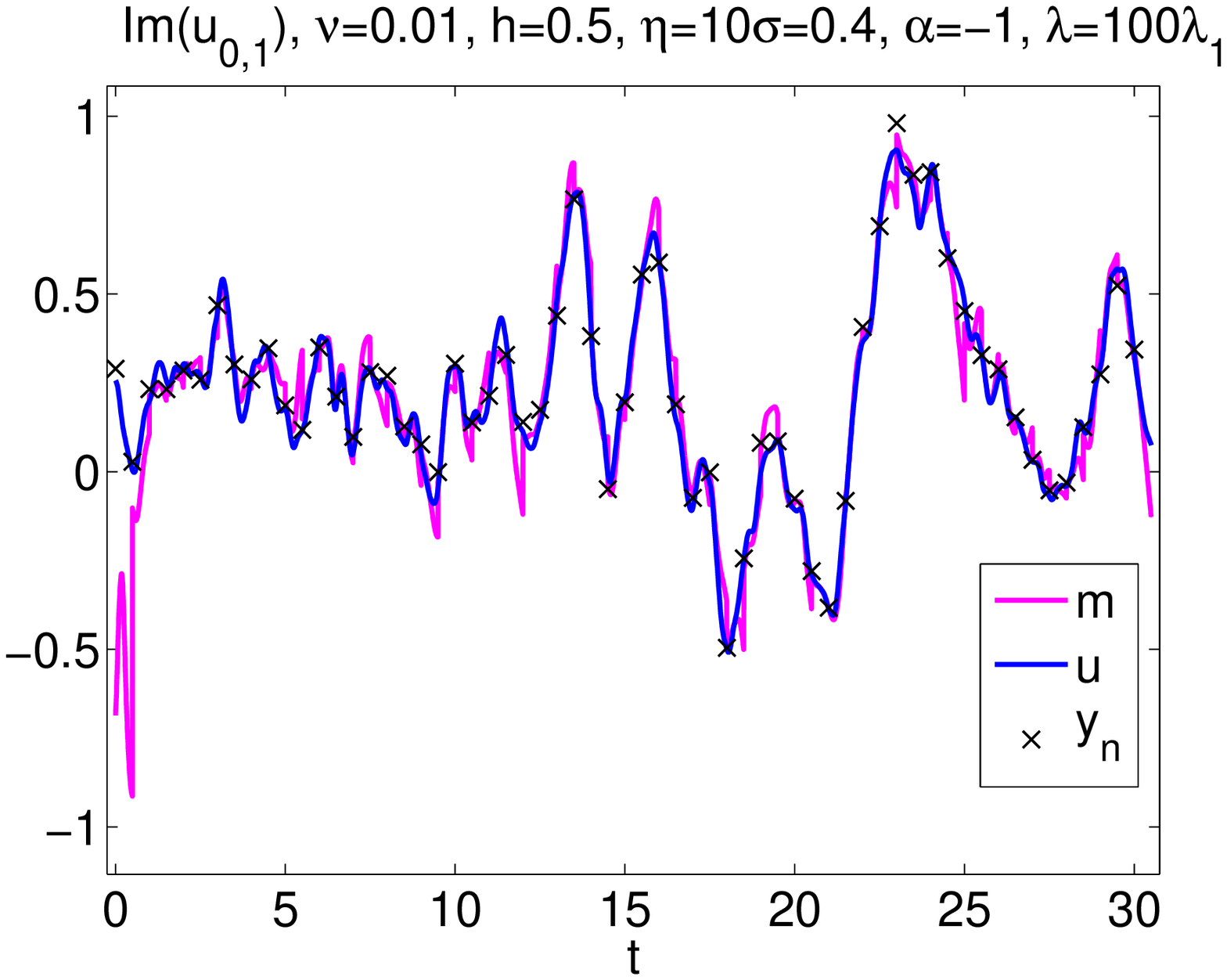}
\includegraphics[width=.45\textwidth]{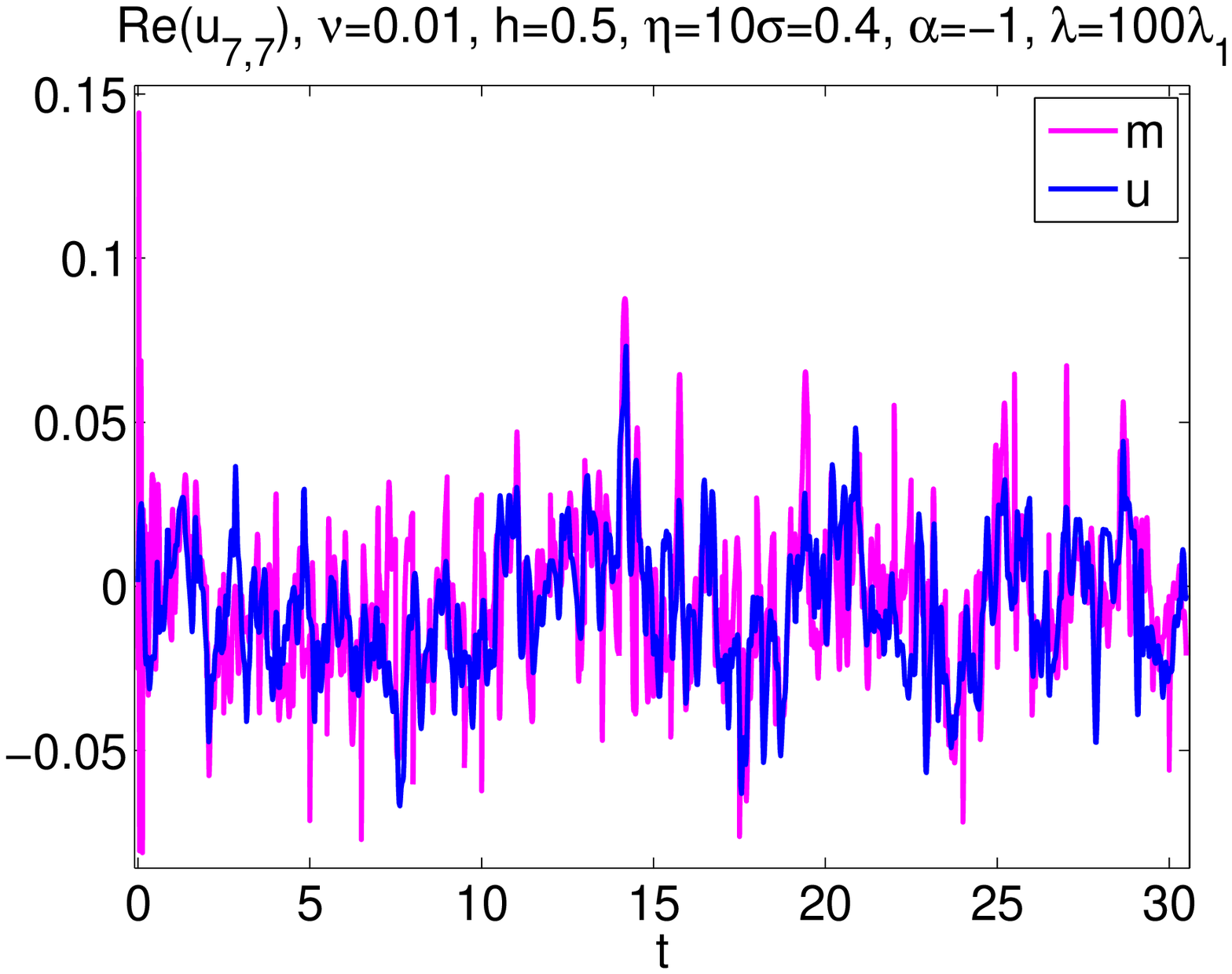}
\includegraphics[width=.45\textwidth]{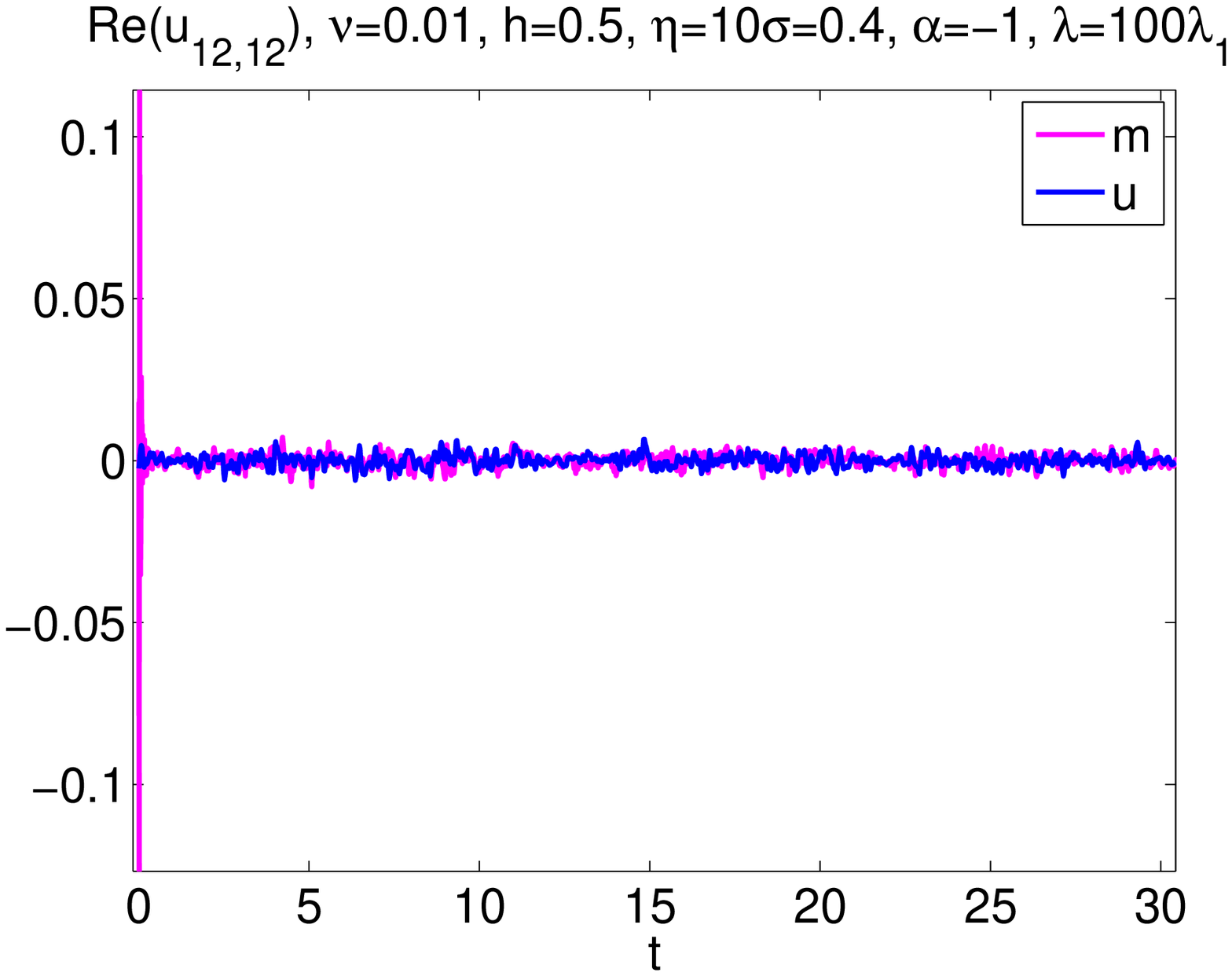}
\caption{Example of an improved estimator for partial observations
  with $\lambda=100 \lambda_1$ and otherwise the  
same parameters as in Fig. \ref{am1.10}. Panels are the same.}
\label{am1.10.p10}
\end{figure*}

We now proceed to examine the case of partial observations,
illustrating Theorem \ref{t:m}.
Note that the forced mode has magnitude $|k_f|^2 =50$,
so ensuring that it is observed requires 
that $\lambda > 50\lambda_1$.
When enough modes are retained, for example when 
$\lambda=100\lambda_1$ in our setting, 
the results for the $\alpha=1$ case remain roughly 
the same and are not shown.  However, in the case $\alpha=-1$, 
in which the observations are trusted more than
the model at high wavevectors, 
the results are greatly improved by ignoring the observations of the
high-frequencies.  See Fig. \ref{am1.10.p10}.  
This improvement, and indeed the improvement beyond setting $B=0$
for both cases $\alpha=\pm 1$ disappears as $\lambda$ is decreased.
In particular, when $\lambda=25 \lambda_1$ the error is never very 
much smaller than the upper bound.  
This is due to the fact that the
dynamics of the low wavevectors tend to be unpredictable and
they contain very little useful information for the assimilation.
Then, for much smaller $\lambda=4 \lambda_1$,
once enough unstable modes are left unobserved,
there is no convergence.  
% The phenomenon is of course improved by considering
% higher frequency in time observations (smaller $h$), 
% since this lessens the manifestation of instabilities.
The order of magnitude of the error in the asymptotic 
regime as a function
of $\eta$ remains roughly consistent as $\lambda$ is decreased
until the estimator no longer converges.  
%, which tends to be close to the
%lower bound. 
For small $h$ (high-frequency in time observations)
and complete observations,
the estimator can be slow to converge to the 
asymptotic regime.
In this case, the number of iterations until convergence, for a sufficiently 
small $\eta$, becomes significantly larger as $\lambda$ is decreased
(again until the estimator fails to converge at all).

\begin{figure*}
\includegraphics[width=.45\textwidth]{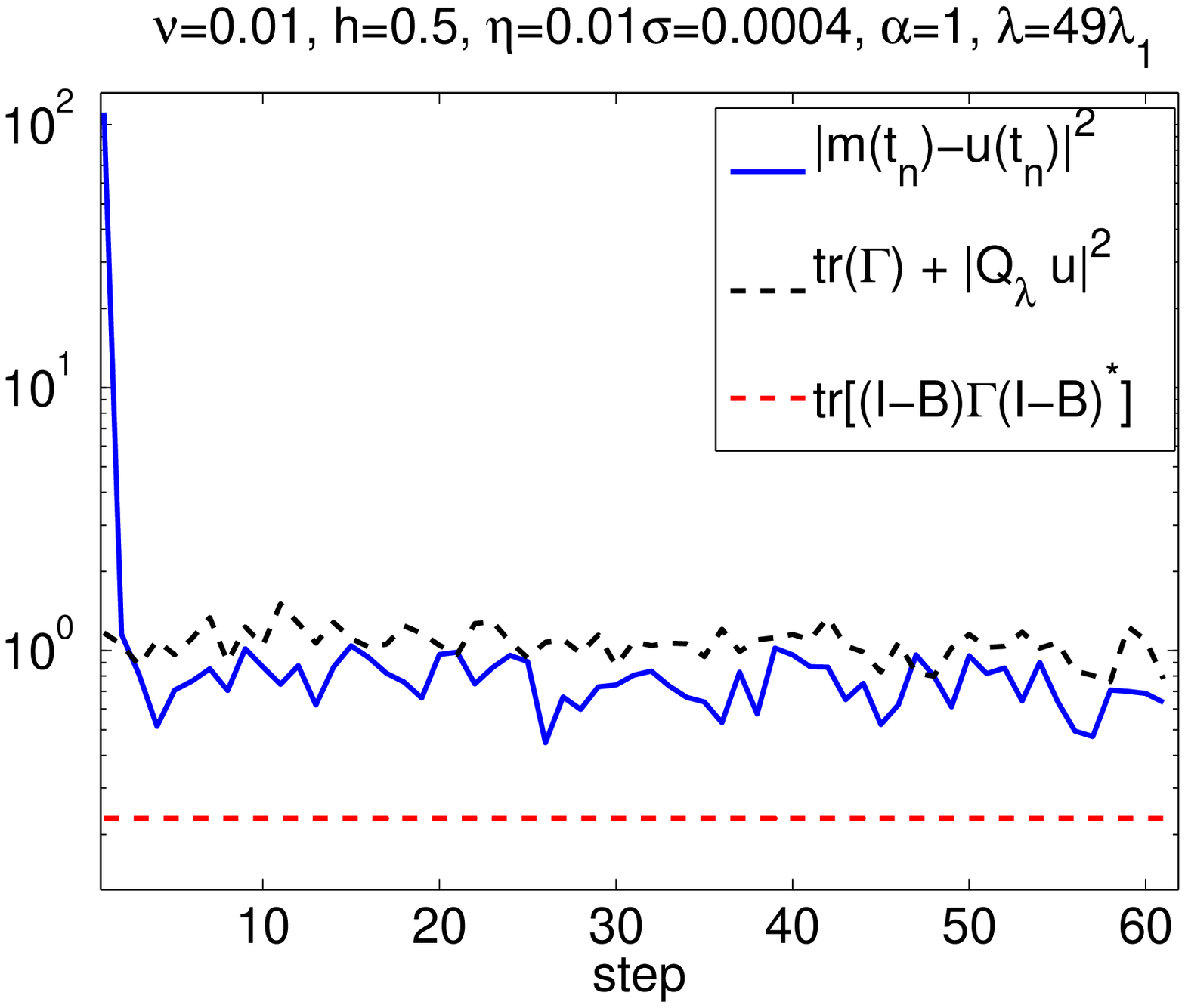}
\includegraphics[width=.45\textwidth]{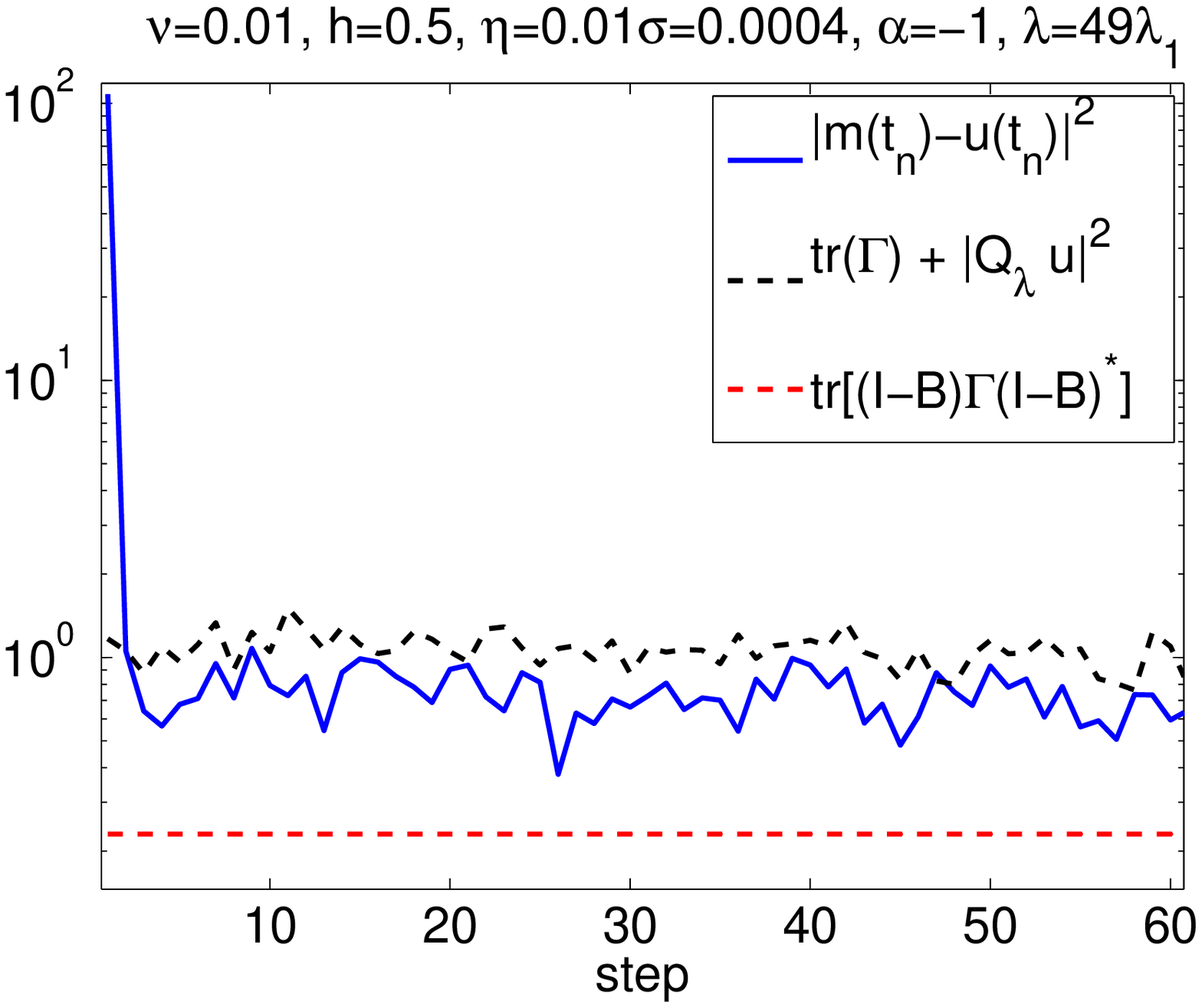}
\includegraphics[width=.45\textwidth]{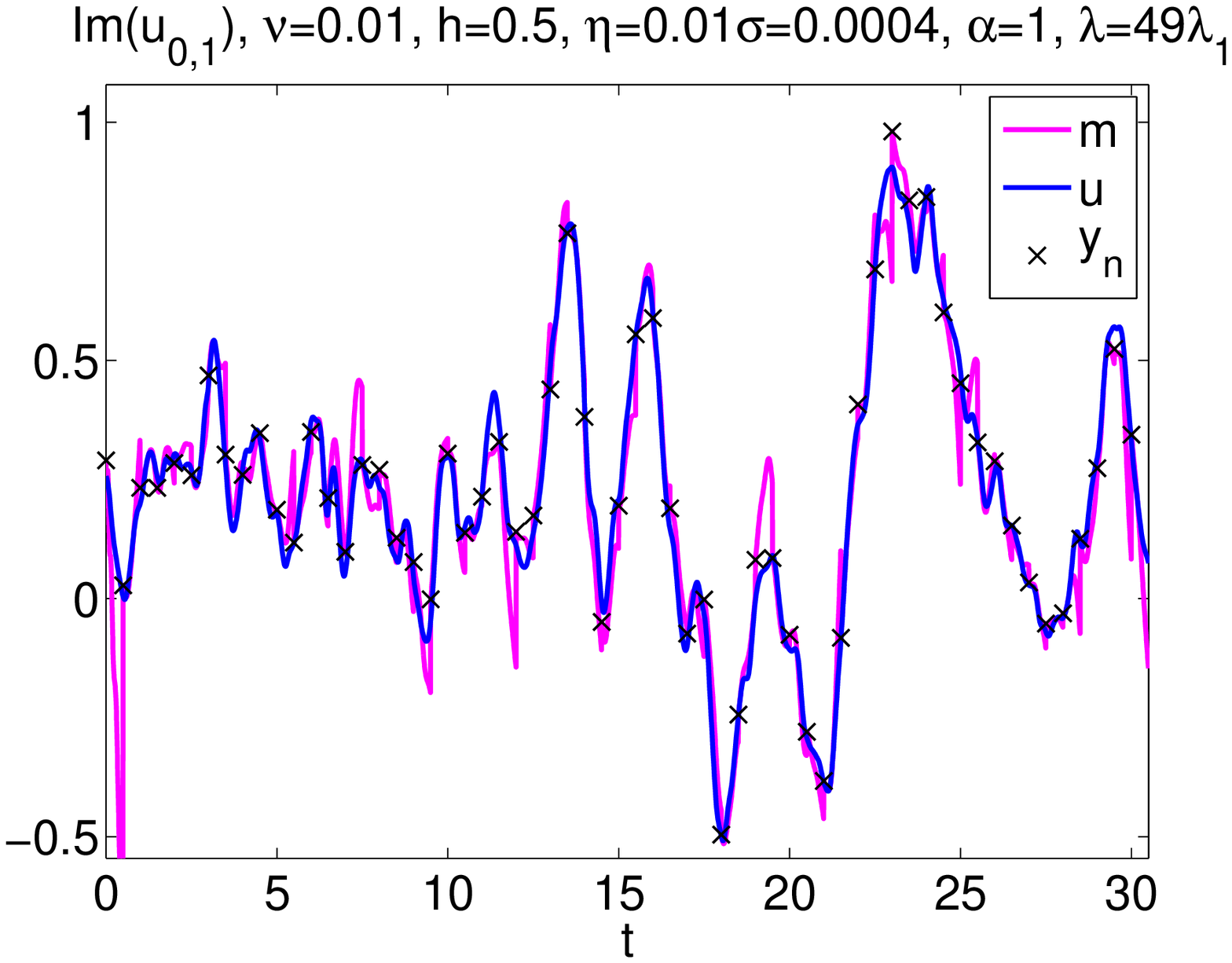}
\includegraphics[width=.45\textwidth]{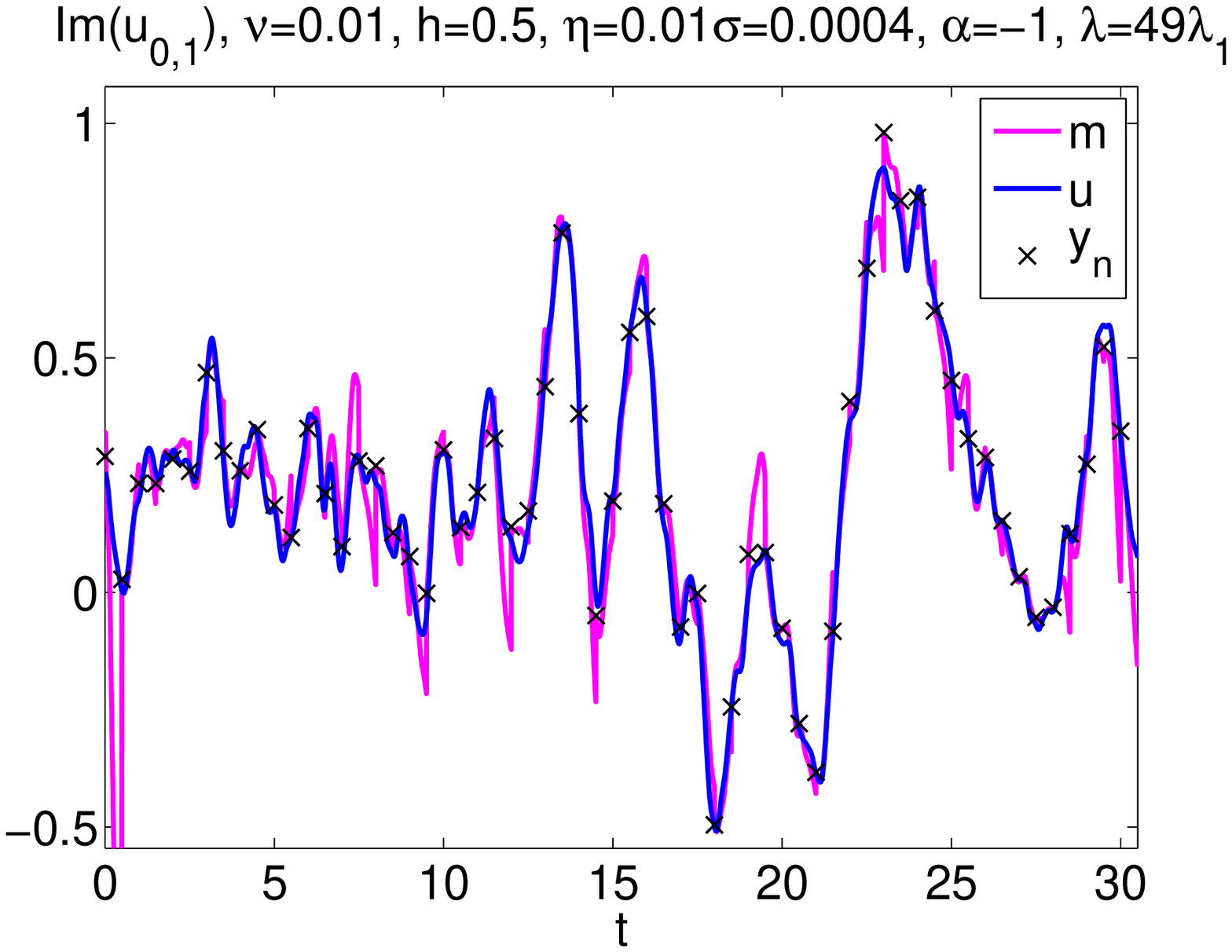}
\caption{Examples of estimators for partial observations 
with $\lambda=49 \lambda_1^2$ and $\eta=0.1\sigma=0.004$, otherwise the  
same parameters as in Figs. \ref{a1.1} and \ref{am1.1}. Left panels
are for $\alpha=1$ and right panels are for $\alpha=-1$.}
\label{a1am1.7.p01}
\end{figure*}

Given $\lambda \approx k_\lambda^2\lambda_1$, we expect that
for $\eta$ sufficiently small the contribution of the model
to the filter will be negligible for all 
$k$ with $|k|<k_\lambda$ for $\alpha=1$.  Hence
the estimators for both $\alpha=\pm 1$ will behave similarly. 
An example of this is shown in Fig. \ref{a1am1.7.p01}
where $\eta=0.01\sigma=0.0004$ 
and $\lambda=49 \lambda_1$ 
in Fig. \ref{a1am1.7.p01}.  In both cases, the estimator is
essentially utilizing all the available observations. 
There are enough observations to draw the 
higher wavevectors of the estimator
closer to the truth than if we just set the population of 
those modes to zero.  
In contrast, as mentioned above, when $\lambda=25 \lambda_1$,
there are not enough observations even when they are all used, 
and the error is roughly the same as the upper bound 
as $\eta \rightarrow 0$ (not shown).

%In the case that too few modes are retained, for example when
%$\lambda = 49 \lambda_1$ in our case, 
%essentially all the observations are needed
%for a reasonable reconstruction and the results for $\alpha=1$ and 
%$\alpha=-1$ are essentially the same: $\eta$ must be 
%sufficiently small that the observations guide the solution in 
%those modes that are retained (the critical value of $\eta$ is 
%slightly larger for the $\alpha=-1$ case).  
%Note that $\lambda < k_f^2\lambda_1$  in this case
%%$k_f \approx \lambda$ 
%so the unstable forcing mode is never observed.

\section{Conclusion}
\label{sec:conclusions}

This paper contains three main components:

\begin{itemize}

\item we show that the filtering problem for the
Navier-Stokes equation may be formulated as a
a sequence of well-posed inverse problems:
Theorem \ref{t:min};

\item we prove Theorems \ref{t:m} and \ref{t:mz},
which establish filter accuracy and stability 
provided variance inflation is employed;

\item we describe 
numerical results which illustrate, and extend the
validity of, the theory.

\end{itemize}

We note that the analysis will also apply to other semilinear dissipative
PDEs which possess a squeezing property (contraction when projected
into the high modes) and a global attractor; such equations
are studied in depth in \cite{book:Temam1997}, 
and include the Cahn-Allen and Cahn-Hilliard equations,
the Kuramoto-Sivashinsky equation and Ginzburg-Landau equation.
These two structural properties of the underlying model, when combined 
with sufficient variance inflation ($\eta$
small enough in 3DVAR) enable proof that the 3DVAR filter has
a contractive property, even when the underlying dynamical system
itself has positive Lyapunov exponents.

There are a number of natural future directions
which stem from this work:

\begin{itemize}

\item to develop analogous filter stability theorems
for more sophisticated filters, such as the extended
and ensemble-based methods, when applied to the
Navier-Stokes equation;

\item to study model-data mismatch by looking at filter
stability for data generated by forward models which
differ from those used to construct the filter; 

\item to study the effect of filtering in the presence of
model error, by similar methods, to understand how this
may be used to overcome problems arising from model-data mismatch.

\item to combine the analysis in this paper, which
concerns nonlinear problems, but assumes that the
observation operator and the Stokes operator commute,
and the recent work \cite{pot12} which concerns only linear
autonomous problems but does not assume commutativity of the
solution operator for the forward model and 
the observation operator.

\end{itemize}

\vspace{0.1in}

\noindent{\bf Acknowledgements}
{AMS would like to thank the following institutions for
financial support: EPSRC, ERC, ESA, and ONR; KJHL was supported
by EPSRC, ESA, and ONR; and 
CEAB, KFL, DSM and MRS were supported EPSRC, through the
MASDOC Graduate Training Centre at Warwick University. 
The authors also thank The Mathematics Institute and Centre for Scientific
Computing at Warwick University for supplying valuable computation
time.  Finally, the authors
thank Masoumeh Dashti for valuable input.}

% \appendix
% \renewcommand{\theappxlem}{A.\arabic{appxlem}}
% \section{Appendix}

%\section*{References}
%\bibliographystyle{iopart-num}
%\bibliographystyle{jphysicsB}
%\bibliographystyle{plain}
%\bibliographystyle{elsarticle-num}
%\biboptions{longnamesfirst,angle,semicolon}

%\bibliography{fsbib}

\end{document}